\documentclass[10pt]{amsart}
\usepackage{color}
\usepackage{soul}
\usepackage{cancel}

\usepackage{amssymb}

\definecolor{c20}{rgb}{0.,0.7,0.}
\definecolor{c30}{rgb}{0.,0.,1.}
\definecolor{c40}{rgb}{1,0.1,0.7}
\definecolor{c50}{rgb}{1,0,0}
\definecolor{c60}{rgb}{1,0.9,0.1}

\newcommand{\kb}[1]{\boldsymbol{#1}}
\newcommand{\vk}[1]{\kb{#1}}

\newcommand{\abs}[1]{\left\lvert #1 \right\rvert}

\newcommand{\E}[1]{\mathbb{E}\left\{#1\right\}}

\newcommand{\pk}[1]{\mathbb{P} \left\{ #1 \right\} }

\newcommand{\R}{\mathbb{R}}

\newcommand{\N}{\mathbb{N}}

\newcommand{\BQN}{\begin{eqnarray}}
\newcommand{\EQN}{\end{eqnarray}}
\newcommand{\BQNY}{\begin{eqnarray*}}
\newcommand{\EQNY}{\end{eqnarray*}}

\newcommand{\BS}{\begin{sat}}
\newcommand{\ES}{\end{sat}}
\newcommand{\BT}{\begin{theo}}
\newcommand{\ET}{\end{theo}}
\newcommand{\BK}{\begin{korr}}
\newcommand{\EK}{\end{korr}}

\newcommand{\BD}{\begin{de}}
\newcommand{\ED}{\end{de}}

\newcommand{\BIT}{\begin{itemize}}
\newcommand{\EIT}{\end{itemize}}
\newcommand{\BDI}{\begin{description}}
\newcommand{\EDI}{\end{description}}

\newcommand{\BRM}{\begin{remarks}}
\newcommand{\ERM}{\end{remarks}}

\newcommand{\BEL}{\begin{lem}}
\newcommand{\EEL}{\end{lem}}

\newtheorem{theo}{Theorem}[section]
\newtheorem{sat}[theo]{Proposition}
\newtheorem{de}[theo]{Definition}
\newtheorem{lem}[theo]{Lemma}

\newtheorem{korr}[theo]{Corollary}
\newtheorem{remark}[theo]{Remark}
\newtheorem{remarks}[theo]{Remarks}
\newtheorem{prop}[theo]{Proposition}

\newcommand{\prooflem}[1]{\textsc{\bf Proof of Lemma} \ref{#1}:}

\newcommand{\COM}[1]{}

\allowdisplaybreaks[1]

\newcommand{\QED}{\hfill $\Box$}

\def\vp{\varepsilon}

\def\rw{\rightarrow}

\def\IF{\infty}

\def\Cov{\mathrm{Cov}}

\date{}

\def\oo{(1+o(1))}

\def\LT{\left}
\def\RT{\right}
\def\H{\mathcal{H}}

\def\rw{\rightarrow}

\def\vn{\varepsilon}
\def\Var{\text{Var}}

\def\todis{\overset{d}\rightarrow}

\def\NN{\mathcal{N}}

\def\ovX{\overline X_H}

\def\bnbH{b_n^{H/\beta}}
\def\bnH0{b_n^{H_0/\beta}}
\def\wY{\widetilde{Y}}
\def\wb{\widetilde{b}}
\def\wa{\widetilde{a}}
\def\wR{\widetilde{R}}
\def\wQ{\widetilde{Q}}
\def\wM{\widetilde{M}}

\usepackage{soul}
\usepackage{color}

\newcommand{\tb}[1]{{\textcolor{blue}{#1}}}

\def\tb#1{#1}

\newcommand{\Ji}[1]{{\textcolor{red}{#1}}}
\newcommand{\Pe}[1]{{\textcolor{blue}{#1}}}

\def\Pe#1{#1}
\def\Ji#1{#1}

\begin{document}

\title[
Extreme value theory for supremum] 
{Extreme value theory for a sequence of suprema of  a class of  Gaussian processes with trend}

 \author{Lanpeng Ji}
\address{Lanpeng Ji, School of Mathematics, University of Leeds, Woodhouse Lane, Leeds LS2 9JT, United Kingdom
}
\email{l.ji@leeds.ac.uk}

\author{Xiaofan Peng}
\address{Xiaofan Peng, School of Mathematical Sciences, University of Electronic Science and Technology of China, Chengdu 611731, China}
\email{xfpeng@uestc.edu.cn}

\bigskip

\date{\today}
 \maketitle

{\bf Abstract:} We investigate extreme value theory of a class of random sequences defined by the all-time suprema of aggregated self-similar Gaussian processes with trend. This study is motivated by its potential applications in various areas and its theoretical interestingness.  We consider both stationary sequences and non-stationary sequences obtained by considering whether the trend functions are identical or not. We show that a sequence of suitably normalised $k$th order statistics converges in distribution to a limiting random variable which can be a negative log transformed Erlang distributed random variable,  a Normal  random variable or a mixture of them, according to three conditions deduced through the model parameters. Remarkably, this phenomenon resembles that for the stationary Normal sequence. We also show that various moments of the normalised $k$th order statistics  converge to the  moments of the corresponding limiting random variable. The obtained results enable us to analyze various properties of these random sequences, which reveals the  interesting particularities of this class of random sequences in extreme value theory.

\medskip

\medskip


\section{Introduction}


Let  $\{X_i(t), t\ge0\}, i=1,2,\ldots,$ be independent copies of a centered self-similar Gaussian process with almost surely (a.s.) continuous sample paths,   self-similarity index $H\in(0,1)$ and variance function $t^{2H}$, and let  $\{X(t), t\ge0\}$ be another independent  centered self-similar Gaussian processes with a.s. continuous sample paths,  self-similarity index $H_0\in(0,1)$ and variance function $t^{2H_0}$.  We refer to the recent contribution \cite{DT20} for a nice discussion on properties and examples of self-similar Gaussian processes.  
We define, for positive constants $\sigma,\sigma_0, c_i, i\ge1$ and $\beta>\max(H,H_0)$,
\BQN \label{eq:Qi}
 Q_i :=\sup_{t\ge 0}( \sigma X_i(t)+\sigma_0 X(t)-c_it^\beta), \ \ \ \ i=1,2,\ldots.
\EQN
Studies on distributional properties of a single  extremum, $Q_i$, or  joint extrema (e.g., $(Q_1,   \ldots,Q_d)$ for some $d\in\N$) have attracted growing interest in recent literature.
On one side, it is a natural object of interest in the extreme value theory of stochastic processes. On the other side, strong motivation for this investigation stems, for example, from (multivariate) stochastic models applied to modern risk theory, advanced communication networks and financial mathematics, to name some of the applied-probability areas. Since explicit formulas for the (joint) distributions are out of reach except for some very special cases, current studies have been focused on deriving  (joint) tail asympotitics for these (joint) extrema; see e.g., \cite{ DHL17,dieker2005extremes,HP99, HS06} on the single extremum and \cite{DJR20,RolskiSPA, KZ16} on the joint extrema of some Gaussian processes.

In this paper, we are interested in developing extreme value theory for the sequence of random variables $\{Q_i\}_{i\ge 1}$ defined in \eqref{eq:Qi}. Precisely, denoting the {\it $k$th largest value}  (or {\it $k$th order statistics}) of $Q_1,   \ldots, Q_n$ by
\BQN \label{eq:MQ}
M_n^{(k)}:=\overset{(k)} {\underset{i\le n}{\max} }\ Q_i,\ \ \ \ 1\le k< n\in \N,
\EQN
we aim to establish limit theorems for $e_n^{-1} ( M_n^{(k)} - d_n)$, as $n\to\IF,$  for some suitably chosen normalising deterministic functions $e_n, d_n, n\in \N$. To this end, it is easily seen that we  can simply assume $\sigma=1$ without loss of generality. For notational simplicity, we will make this assumption throughout the rest of the paper but we should bare in mind that general results for $\sigma\neq 1$ can be derived by an easy adjustment of the normalizing functions $e_n, d_n$. Note that by our notation  $M^{(1)}_n=\overset{(1)} {\underset{i\le n}{\max} }\ Q_i$ should be understood as $M_n=\max_{i\le n}Q_i$. 

One of motivations for this study stems from the increasing interest in developing extreme value theory for Gaussian processes in random environment. We refer to a series of papers by Piterbarg and his co-authors, e.g., \cite{MR2832415}, \cite{Pit12} and \cite{PPS17}, for recent developments on this topic.  The model in \eqref{eq:Qi} can be seen as a multi-dimensional counterpart and may be  interpreted as follows: The processes $\{\sigma X_i(t)-c_it^\beta, t\ge 0\}$, $i=1,2,\ldots, n,$ are the subject of primary interest that are affected (or perturbed)  by a common random environment modelled by $\{\sigma_0X(t), t\ge 0\}$. For example, we consider an insurance company running $n$ lines of business, where the  claim surplus processes are modelled by $\{\sigma X_i(t)-c_it^\beta, t\ge 0\}$, $i=1,2,\ldots, n$, which are affected by the common  random environment  $\{\sigma_0 X(t), t\ge 0\}$.
In this context, one could think of  environmental factors
affecting the claim surplus processes, such as the state of the economy, the political situation,
weather conditions, and policy regulations. For this model, 
the obtained results below can give an approximation for the probability
$$
\pk{M_n^{(k)} > d_n + e_n x}
$$
for any $x\in \R$, and large $n$. This approximation can be used to evaluate the probability that at least $k$ lines of business will ultimately get ruined, if ruin is defined to occur when the claim surplus process exceeds the capital $d_n + e_n x$ for a chosen $x$ and a chosen large $n$. 

 Another motivation of our study comes from a recent contribution \cite{MSJVZ21}, where the authors investigate a particular model with all the Gaussian processes involved in \eqref{eq:Qi}-\eqref{eq:MQ} being Brownian motions, $c=c_i, i\ge 1,$ and $\beta=1$ (hereafter, referred to as {\it Brownian model with linear drift}).  In their context, $M_n$ 
models the maximum queue length in
a fork-join network of $n$ statistically identical queues (i.e., $Q_1,  \ldots, Q_n$) which are driven by a common Brownian motion perturbed arrival process and  independent  Brownian motion perturbed service processes, respectively. The obtained theoretical limit result on $M_n$ therein is the key to developing structural insights into the dimensioning of assembly systems; interested readers are referred to \cite{MSJVZ21} for more details on this application. We extend the model of  \cite{MSJVZ21} by considering  general self-similar Gaussian processes with (non-identical) non-linear drifts, and study the limit properties for the $k$th order statistics of $Q_1,\ldots, Q_n$. It is known that Gaussian processes play an important role in network modelling, see e.g., \cite{Man07}, we expect that the obtained results  potentially have general applications in this area, 
particularly 
for fork-joint networks.


Besides the motivations in the application areas as above, the study of the random sequence $\{Q_i\}_{i\ge 1}$ is of  interest from a theoretical point of view. Note that the random variables  $Q_i, i\ge 1$ are mutually dependent with a  dependence structure induced by the common stochastic process  $\{\sigma_0 X(t), t\ge 0\}$.   The random sequence $\{Q_i\}_{i\ge 1}$ is a {\it stationary  sequence} if $c=c_i, i\ge 1$, and a {\it non-stationary sequence}, otherwise. 
Most of the work in extreme value theory has been done for independent and identically distributed  (IID)  sequences; see, e.g., \cite{DeHF06,EKM97,Res87}.  
We refer to Chapters 3 and 5-6 in \cite{leadbetter1983extremes} and Chapters 8-9 in \cite{FHR10} for some early discussions on extreme value theory for general dependent (stationary or non-stationary)  sequence. Extreme value theory for general class of (non-)stationary  sequences 
normally involves an asymptotic independence condition $D(u_n)$ of mixing type. Clearly, the stationary sequence $\{Q_i\}_{i\ge 1}$, with $c=c_i, i\ge 1$, is non-ergodic (in fact, exchangable) and thus the mixing condition $D(u_n)$ can not be verified. This means that many results in the classical theory may not be proved using the classical point process approach. For example, a Poisson point process approximation for the number of high-level exceedances by $\{Q_i\}_{i\ge 1}$ may be impossibe (or meaningless as commented in Example 8.2.1 in \cite{FHR10}), since a proof of such a result normally relies on the mixing condition.  

Although the dependence structure of the sequence  $\{Q_i\}_{i\ge 1}$ is generally hard to analyze (see \cite{DJR20} for some remarks on the Brownian model with linear drift), this sequence exhibits some interesting (sometimes uncommon) properties as follows:

\begin{itemize}

\item  
Some convergence results for suitably normalised order statistics $e_n^{-1} ( M_n^{(k)} - d_n)$, as $n\to\IF$, are approachable directly without first deriving point process convergence as in the classical theory.  
More precisely, we show, in Theorem \ref{Thm:1} and Theorem \ref{Thm:1-1} that, for different scenarios (derived trough $H,H_0,\beta$), there exist $e_n,d_n,n\in \N$ such that the weak limit of $e_n^{-1} ( M_n^{(k)} - d_n)$ can be a negative log transformed Erlang distributed random variable,  a Normal  random variable or a mixture of them. This phenomenon remarkably resembles that of the stationary Normal sequences, for which it has been well known that  a sequence of suitably normalised maxima of a stationary Normal sequence with correlation function $r_n$ will converge to a Gumbel random variable, a Normal random variable or a  mixture of them, according to different limiting values of $r_n \log n$; see \cite{leadbetter1983extremes}. By considering the $k$th order statistics we can obtain an equivalent (mixed) Poisson distribution convergence of the number of high-level exceedances by $\{Q_i\}_{i\ge 1}$; some result that might be useful in applications. 

\item The notion of phantom distribution function was introduced by O'Brien \cite{OB87}, and that the existence of such a distribution is a quite common phenomenon for staitonary weakly dependent sequences.  In \cite{DJL15} the authors derive equivalence statements for the existence of a continuous phantom distribution function for a stationary sequence, where it is claimed that the asymptotic independence of maxima  (i.e., $D(u_n)$) is not really a necessary condition. They also constructed a non-ergodic stationary process which admits a continuous phantom distribution function; see Theorem 4 therein. The stationary sequence $\{Q_i\}_{i\ge 1}$ (with $c=c_i, i\ge 1$) here gives another example of non-ergodic stationary process which admits a continuous phantom distribution function under some scenario, while under other scenarios it does not admit a  phantom distribution function. See Remarks \ref{Rem:1} (d) below for some detailed discussions.

\item  It is well known that for IID sequences, the weak convergence of any one of the normalised order statistics is equivalent to the convergence of the corresponding maxima; see e.g., Theorem 2.2.2 in \cite{leadbetter1983extremes}. However, for general dependent random sequences this might not always be true, for example, Mori \cite{Mo76}   provides an example of such a sequence for which convergence for lower order statistics does not guarantee the convergence for the higher order statistics.  We refer to \cite{Hs88}  for some related discussions. 
In particular, 
the author  derives an interesting equivalence result on the (compound) Poisson point process approximation for exceedences and the convergence for all the order statistics (see Theorem 5.1 therein). Note again that a mixing condition is  crucial in the proof of this result. 
The sequence $\{Q_i\}_{i\ge 1}$ here, regardless of stationarity, gives an example where 
the convergence of all the order statistics can be established. Since the results obtained in  Theorem \ref{Thm:1} and Theorem \ref{Thm:1-1} are scenario specific, it is easily seen that the convergence for  any one of the order statistics  is equivalent to the convergence for the maxima. 

\item An overview of extreme value theory for general non-stationary sequences can be found in Chapter 9 in \cite{FHR10} which was developed on the basis of assuming some general mixing condition $D(u_{ni})$. As discussed earlier, the non-stationary sequence $\{Q_i\}_{i\ge 1}$ (with different $c_i$'s) may not satisfy this mixing condition.  That being said, we can  derive, under a mild  restriction on $c_i$'s, some convergence results for the order statistics $ M_n^{(k)}$ with some suitable normalisation, which can be seen as a thinning version of the results for stationary cases where all the $c_i$'s are the same.

\end{itemize}

\COM{
As an exceptional and tractable example, stationary Normal sequence, say $\{Y_i\}_{ i\ge 1},$  was discussed in more detail in the literature;
see, e.g.,  \cite{MY75} or Chapter 6 in \cite{leadbetter1983extremes}. Denote  by $r_n=\text{Cov}(Y_1, Y_{n+1}), n\ge 1,$ the covariance function of $\{Y_i\}_{ i\ge 1}$. It has been
shown 
that  if the weak dependence condition  $\lim_{n\to\IF} r_n \log n =0$ (also called Berman condition)  holds, then a sequence of suitably normalised  $\max_{i\le n} Y_i$  converges to a Gumbel random variable. Whereas, if $\lim_{n\to\IF} r_n \log n = \IF$ holds, then a   sequence of suitably normalised  $\max_{i\le n} Y_i$   converges  to a Normal random variable. Moreover,
  if the strong dependence condition $\lim_{n\to\IF} r_n \log n =r \in(0,\IF)$ holds, then  a sequence of suitably normalised  $\max_{i\le n} Y_i$  converges   to a mixture of a Gumbel random variable and a Normal random variable.  This reveals that the class of limiting random variables of the normalised maxima for dependent random variables is much richer. Moreover, we can see from \cite{leadbetter1983extremes} that some weak convergence results were also established for suitably normalised  order statistics of $\{Y_i\}_{i\ge 1}$, under  weak or strong dependence conditions.
}

It is well known that \Pe{weak} convergence of a sequence of random variables does not imply  convergence of moments. In the classical extreme value theory for IID sequences, it has been shown that such convergence of moments of normalised maxima is valid provided that some moment conditions are satisfied; see, e.g., \cite{Pick68} or Section 2.1 of \cite{Res87}. More recently, in  \cite{BO13} and \cite{PNS13} the authors discuss  moment convergence of extremes under power normalization. A natural question here is whether  various moments of $e_n^{-1} ( M_n^{(k)} - d_n)$ converge to the  moments of the corresponding limiting random variable. In Theorem \ref{Thm:2} (see also Theorem \ref{Thm:1-1}), this question is answered affirmatively without imposing any further conditions. As a by-product of this study, we derive, in Proposition \ref{Lem-mom-Ynk}, some moment convergence result for the $k$th order statistics of  IID generalized Weibull-like random variables. This  problem is of independent interest and  has not  been extensively explored in the existing literature as far as we are aware.

It is worth mentioning that in \cite{MSJVZ21}, the stationary and independent increment property of Brownian motion is the key to their proof for the weak convergence.  Whereas, our proof mainly relies on the asymptotic theory of Gaussian processes, particularly, the Borell-TIS inequality and the tail asymptotics of supremum. Some ideas in our proof for the weak convergence results 
are stimulated by the intuitive interpretations 
provided for the Brownian model with linear drift in Section 5.1 of \cite{MSJVZ21}.

\medskip

Brief outline of the paper: In Section 2 we  present some preliminary results concerning the tail asymptotics of the all-time supremum of  a class of self-similar Gaussian processes with trend and some limit results for  generalized Weibull-like random variables. The main results on the convergence of the suitably normalised order statistics are given in Section 3, with the proofs displayed in Section 4. Some technical proofs for Section \ref{Sec:Gam} and some frequently used $C_r$ inequalities are presented in an Appendix.

\section{Preliminaries}


\subsection{Extremes of self-similar Gaussian processes with trend}\label{s.notation}

Let $\{X_H(t), t\ge0\}$ be  a centered self-similar Gaussian process with a.s. continuous sample paths, self-similarity index $H\in(0,1)$ and variance function $t^{2H}$. 
Define, for $\beta>H$ and $c>0$,
$$
\tau_u \ =  \ \inf\{t\ge0 : X_H(t) -ct^\beta >u\}
$$
to be the first hitting time of a level $u>0$  by the stochastic process $\{X_H(t) -ct^\beta, t\ge 0\}$.

By self-similarity of $X_H$,  the probability of ultimately crossing an upper level $u>0$ by the process $\{X_H(t) -ct^\beta, t\ge 0\}$  is given as
\BQNY 
\pk{\tau_u<\IF}=\pk{\sup_{t\ge0} \Bigl(X_H(t)-ct^\beta\Bigr)>u}=
\pk{\sup_{t\ge0}  Z(t) >u^{1-{H}/{\beta}}},
\EQNY
with
\BQNY
Z(t)=\frac{X_H(t)}{1+ct^\beta},\ \  t\ge0. 
\EQNY
It follows from Proposition 3 of \cite{dieker2005extremes} that
\BQNY
\lim_{t\to\IF} Z(t) =0\ \ \ a.s.
\EQNY
which means that the sample paths of $\{Z(t),t\ge0\}$ are  bounded   a.s. This ensures  $\sup_{t\ge0} \LT(X_H(t)-ct^\beta\RT)<\IF$ a.s. and is important when we apply the Borell-TIS inequality later in some of the proofs.
Next,  from  \cite{HP99} or \cite{ HP08}, we have that the standard deviation function
$$
\sigma_Z(t)=\sqrt{\Var(Z(t))}=\frac{t^H}{1+ct^\beta}, \ \ \ \ t\ge 0,
$$
attains its maximum on $[0,\IF)$ at the unique point
$t_0=\LT(\frac{H}{c(\beta-H)}\RT)^{\frac{1}{\beta}}$
and
\[
\sigma_Z(t)=A-\frac{BA^2}{2}(t-t_0)^2 +o((t-t_0)^2), \ \ \ \  t\to t_0,
\]
 where
\BQN\label{eq:AB}
 A=\frac{t_0^H}{1+ct_0^\beta}=\frac{\beta-H}{\beta} \LT(\frac{H}{c(\beta-H)}\RT)^{\frac{H}{\beta}},\ \ \ \ B=\LT(\frac{H}{c(\beta-H)}\RT)^{-\frac{H+2}{\beta}} H\beta.
\EQN
Furthermore, we assume a {\it local stationarity}
of the standardized  Gaussian process $\ovX(t):=X_H(t)/t^H, t> 0$ in a
neighbourhood of the point $t_0$, i.e.,
\BQN\label{eq:locstaXH}
\lim_{s,t\to t_0}\frac{ \E{(\ovX(s)-\ovX(t))^2}}{K^2 (\abs{s-t})}=1 
\EQN
holds for some positive function $K(\cdot)$ which is regularly varying at 0 with index $\alpha/2\in(0,1)$.
Condition \eqref{eq:locstaXH} is a common assumption in the literature; see, e.g., \cite{dieker2005extremes} and \cite{HP99}.
It is worth noting that the assumption \eqref{eq:locstaXH} is slightly general than the S2 in the definition of self-similar Gaussian processes in \cite{DT20}, and in \cite{HP99} a slightly larger class of Gaussian processes is also discussed. 
Throughout this paper, we denote by $ \overset{\leftarrow}{K}(\cdot)$  the asymptotic inverse of  $K(\cdot)$, and thus 
 $$
 \overset{\leftarrow}{K}(K(t))=K( \overset{\leftarrow}{K}(t))\oo=t\oo,\ \ \ t\downarrow0.
  $$
It follows that  $ \overset{\leftarrow}{K}(\cdot)$ is regularly varying at 0 with index ${2}/{\alpha}$; see, e.g., \cite{EKM97}.

Below, by $\{B_{\alpha/2}(t), t\ge 0\}$ we denote a standard fractional Brownian motion
 (sfBm)
with Hurst index $\alpha/2\in(0,1)$, and 
\BQNY
\Cov(B_{\alpha/2}(t),B_{\alpha/2}(s))=\frac{1}{2}(t^{\alpha}+s^{\alpha}-\mid t-s\mid^{\alpha}),\ \ \ t,s\ge0.
\EQNY
The well known Pickands constant $\H_{\alpha}$ in the Gaussian theory is defined by
\[
\H_{\alpha}=\lim_{T\to\infty}\frac{1}{T}
\E{\exp \left(\sup_{t\in[0,T]}(\sqrt{2}B_{\alpha/2}(t)-t^\alpha)\right)} \in (0,\IF).
\]
We refer to  \cite{DH20, demiro2003simulation,  DikerY, Pit96} and references therein
for basic properties of {the} Pickands and related constants.

The following proposition gathers some  useful results from  \cite{HP99} and \cite{HP08} (see also \cite{DHJ15}). 

\begin{prop} \label{prop1} Let $\{X_H(t), t\ge0\}$ be  a centered self-similar Gaussian process defined as above satisfying \eqref{eq:locstaXH}  and let $c>0$. Assume $\beta>H$. Then, for  any $\vn_0\in(0,t_0)$ and
any $T>t_0,$
\BQN\label{eq:HPA1}
\pk{\sup_{t\ge 0}   \Bigl(X_H(t)-ct^\beta\Bigr)>u}&=&\pk{\sup_{t_0-\vn_0\le t \le t_0+\vn_0}   \frac{X_H(t)}{1+ct^\beta}>u^{1-H/\beta}}\oo \nonumber\\
&=&\pk{\sup_{0\le t \le T}   \frac{X_H(t)}{1+ct^\beta}>u^{1-H/\beta}}\oo\nonumber\\
&=&R(u)\exp\LT(-\frac{u^{2\LT(1-\frac{H}{\beta}\RT)}}{2A^2}\RT)\oo, \ \ \ \  u\to\IF,
\EQN
where (with $A,B$ given  in \eqref{eq:AB}) 
\BQN\label{Ru}
R(u)=\frac{A^{\frac{3}{2}-\frac{2}{\alpha}}  \H_\alpha}{2^{ \frac{1}{\alpha}} B^{\frac{1}{2}}}\frac{ u^{  \frac{2H}{\beta}-2}}{ \overset{\leftarrow}K (u^{\frac{H}{\beta}-1}) }, \ \  u>0.
\EQN
\COM{Furthermore, as $u\rw\IF$
\BQN\label{eq:tau}
\frac{\tau_u-t_0 u^{\frac{1}{\beta}}}{A^{\frac{1}{2}} B^{-\frac{1}{2}} u^{\frac{H}{\beta}+\frac1\beta-1}}\Big| (\tau_u<\IF) \todis \NN,
\EQN
with $\NN$ \tb{a standard Normal} random variable.
}
\end{prop}


\subsection{Limit theorems for order statistics of Weibull-like random variables} \label{Sec:Gam}
As in \cite{DEJ14}, a probability distribution function $F$ is called a {\it generalized Weibull-like distribution} if
\BQN\label{eq:Weibull}
F(x) = 1- \rho(x) \exp(-C x^\tau), \ \  \ x\ge x_0
\EQN
for some $x_0>0,$ where $C, \tau$ are two positive constants and  $\rho(x)>0$ is a regularly varying function at infinity with index $\gamma\in \R$.  Note that in \cite{AHLT18, EKM97, GLU15}, the special case $\rho(x)= \rho_0x^\gamma$, for some $\rho_0> 0$,  is discussed.


Let $\{Y_i\}_{i\ge 1}$ be a sequence of IID random variables  which are right tail equivalent to a generalized Weibull-like distribution function of the form
\eqref{eq:Weibull}.
The following result gives a limit theorem for the $k$th order statistics $Y_n^{(k)}:=\overset{(k)} {\underset{i\le n}{\max} }\ Y_i$. 
Hereafter,  $\overset{d}\to$  denotes convergence in distribution and $\overset{d} =$ means equivalence in (finite-dimensional) distribution.



\begin{prop} \label{Prop:GamLim} Let
\BQN\label{def-mun-nun}
\mu_n& =& (C^{-1}\log n)^{1/\tau} +\frac{1}{\tau} (C^{-1}\log n)^{1/\tau-1} \LT( C^{-1} \log(\rho((C^{-1}\log n)^{1/\tau}))\RT), \ n\in \N,\\ \nonumber
\nu_n& =& (C\tau)^{-1} (C^{-1}\log n)^{1/\tau-1}, \ n\in \N.
\EQN
We have, for any fixed integer $k>0,$
\BQNY 
\nu_n^{-1} \LT(Y_n^{(k)} -\mu_n\RT) \overset{d}\to \Lambda^{(k)}, \ \ \ \ n\to\IF,
\EQNY
where  $\Lambda^{(k)} = -\ln E_k$,  with $E_k$ being an Erlang distributed random variable with  shape parameter $k$ and rate parameter equal to 1. In particular, $\Lambda^{(1)}$ is the standard Gumbel random variable, i.e., $\pk{\Lambda^{(1)}\le x} =\exp(-e^{-x}), x\in \R$.

\COM{Equivalently,  for any $x\in\R, n\in \N$,
denote $u_n(x) = \mu_n+\nu_n x$ to be a level and $S_n(x)$ to be the number of exceedances of level $u_n(x)$ by $Y_1, \ldots, Y_n$. We have, for any $k\in \N,$
\BQNY
\lim_{n\to\IF}\pk{S_n(x)<k} = \exp\LT(-e^{-x}\RT) \sum_{l=0}^{k-1}\frac{e^{-lx}}{l !},
\EQNY
that is, the number of exceedances $S_n(x)$ is approximately Poisson distributed with intensity $\lambda(x) = e^{-x}$.
}

\end{prop}



The next result is about the (absolute) moment convergence of the normalized $k$th  order-statistics $\nu_n^{-1} (Y_n^{(k)} -\mu_n)$ defined in Proposition \ref{Prop:GamLim}. To this end, we need to control the left tail of the
generalized Weibull-like random variables $Y_i$.
This problem does not seem to have been explored in the existing literature. 
Some results exist only when $k=1$, that is, for the maximum; see, e.g., \cite{Pick68} or \cite{Res87}.

\begin{prop}\label{Lem-mom-Ynk}
Suppose
\BQN\label{assu-Y-leta}
\limsup_{ x\to\IF}\pk{Y_1<-x}   x^{\eta} <\IF
\EQN
holds for some $\eta>0$. 
We have, for any $\lambda>0$, that
\BQNY
\lim_{n\to\IF}\E{\abs{\nu_n^{-1} \LT(Y_n^{(k)}-\mu_n\RT) }^{\lambda}} = \E{\abs{\Lambda^{(k)}}^{\lambda}}.
\EQNY
\end{prop}


As an application of the above results, we consider order statistics of independent random variables obtained by removing the process $\{\sigma_0 X(t), t\ge0\}$ from \eqref{eq:Qi}-\eqref{eq:MQ} which are defined as
\BQN\label{def:MteQte} 
\wM_n^{(k)}:= \overset{(k)} {\underset{i\le n}{\max} }\  \wQ_i := \overset{(k)} {\underset{i\le n}{\max} }\ \sup_{t\ge 0}(X_i (t)   -c_i t^\beta),\quad  n>  k.
\EQN
Recall that we have assumed  $\sigma=1$.
Without loss of generality, for any fixed $n$ we assume that the constants $c_i$'s are of ascending order with 
\BQN\label{eq:cc}
c:=c_1 =\cdots =c_{m_n}<c_{m_n+1}\le \cdots\le c_n,
\EQN 
where $m_n \le n$ is some integer such that $\lim_{n\to\IF}m_n/n = p\in(0,1]$, i.e., the number of minimal drifts is proportional to the total number $n$. In what follows, 
when we say $m_n =n$ we simply mean that all the $c_i$'s are equal to $c$ and thus  assuming $\{\wQ_i\}_{i\ge 1}$ is an IID  sequence.

Comparing \eqref{eq:HPA1} and \eqref{eq:Weibull}, we see that each $\wQ_i$ is right tail equivalent to a generalized Weibull-like distribution. Particularly, for $\wQ_1$ we have
\BQN \label{eq:tauC}
\rho(u)= R(u),\ \ \ \tau = 2(1-H/\beta), \ \ \ C=\frac{1}{2A^2}.
\EQN
With $A$ given in \eqref{eq:AB}, $R(u)$ given in \eqref{Ru} and $\tau$ given in \eqref{eq:tauC}, we define
\BQN\label{eq:bnan}
b_n& :=& (2A^2 \log n)^{1/\tau} +\frac{1}{\tau} ( 2 A^2 \log n)^{1/\tau-1} \LT( 2A^2 \log(R((2A^2\log n)^{1/\tau}))\RT),\  n\in \N,\\
a_n& :=& 2A^2 \tau^{-1} (2A^2\log n)^{1/\tau-1}, \ n\in \N.   \nonumber
\EQN


\begin{prop} \label{propQQ} Assume that \eqref{eq:cc} holds with some $m_n \le n$ such that $\lim_{n\to\IF}m_n/n = p\in(0,1]$. We have, 
\BQN\label{eq:wMb}
a_{m_n}^{-1}(\wM_n^{(k)} -b_{m_n} ) \overset{d}\to \Lambda^{(k)},  \ \ n\to\IF,
\EQN
and, for any $\lambda>0$,
\BQN\label{eq:wMb2}
\lim_{n\to\IF}\E{\abs{a_{m_n}^{-1}(\wM_n^{(k)} -b_{m_n} ) }^{\lambda}} = \E{\abs{\Lambda^{(k)}}^{\lambda}}.
\EQN
\end{prop}

\section{Main results}

\COM{
We first introduce some additional  notation. Define
\BQNY\label{eq:wMQ}
\wQ_i:=\sup_{t\ge 0}   \Bigl(X_i(t)-ct^\beta\Bigr), \ \ i\in\mathbb N.
\EQNY
Comparing \eqref{eq:HPA1} and \eqref{eq:Weibull}, we see that $\wQ_i$ is right tail equivalent to a generalized Weibull-like distribution with
\BQN \label{eq:tauC}
\rho(u)= R(u),\ \ \ \tau = 2(1-H/\beta), \ \ \ C=\frac{1}{2A^2}.
\EQN
With $A$ given in \eqref{eq:AB},  $R(u)$ given in \eqref{Ru} and $\tau$ given in \eqref{eq:tauC}, we define
\BQN\label{eq:bnan}
b_n& :=& (2A^2 \log n)^{1/\tau} +\frac{1}{\tau} ( 2 A^2 \log n)^{1/\tau-1} \LT( 2A^2 \log(R((2A^2\log n)^{1/\tau}))\RT),\  n\in \N,\\
a_n& :=& 2A^2 \tau^{-1} (2A^2\log n)^{1/\tau-1}, \ n\in \N.   \nonumber
\EQN
}

In this section, we shall first consider the stationary sequence $\{Q_i\}_{i\ge1}$ where $c=c_i, i\ge 1$, and then present results for the general non-stationary case where the constants $c_i$'s may not be the same. Finally, as an application a fractional Brownian model with linear drift is discussed.

Below is one of the principal results on the weak convergence of suitably normalised $M_n^{(k)}$ defined in \eqref{eq:MQ} for the stationary sequence $\{Q_i\}_{i\ge1}$. 
This result extends one of the main results in \cite{MSJVZ21}  where only the  Brownian model with linear drift is discussed.
We also present an equivalent (mixed) Poisson distribution convergence result on the number of exceedances of a level $u_n(x)=b_n+a_n x$  by $Q_1,\ldots, Q_n$, denoted by $N_n(x)$, for any $x\in \R$.
In what follows, we denote by $\NN$ a standard Normal random variable, \Ji{independent of $\Lambda^{(k)}$.}

\BT \label{Thm:1}
Let $M_n^{(k)}, n\in \mathbb N$ be defined in  \eqref{eq:Qi}-\eqref{eq:MQ}  with $\sigma=1$ and    $c=c_i, i\ge1$,  and let $b_n, a_n, n\in \mathbb N$ be given as in
\eqref{eq:bnan}.
 Assume $\beta>\max(H,H_0)$.
We have, for any $k\in \N,$
\begin{itemize}
\item[(i).] If $\beta>2H-H_0$,  then
\BQNY
\sigma_0^{-1}t_0^{-H_0}  b_n^{- H_0 /\beta} (M_n^{(k)} -  b_n)    \ \overset{d}\to \  \NN, \ \ \ \ n\to\IF.
\EQNY

\item[(ii).] If   $\beta<2H-H_0$, then
\BQNY
a_n^{-1} (M_n^{(k)} -  b_n)    \ \overset{d}\to \  \Lambda^{(k)}, \ \ \ \ n\to\IF,
\EQNY
or equivalently, for any $x\in\R,$
\BQNY
\lim_{n\to\IF}\pk{N_n(x)<k} = \exp\LT(-e^{-x}\RT) \sum_{l=0}^{k-1}\frac{e^{-lx}}{l !}.
\EQNY
That is, the number of exceedances $N_n(x)$ is approximately Poisson distributed with intensity $\lambda(x) = e^{-x}$.
\item[(iii).] If   $\beta=2H-H_0$, then

\BQNY
a_n^{-1} (M_n^{(k)} -  b_n)  \ \overset{d}\to \  \Lambda^{(k)} + \frac{\sigma_0c\beta}{H} \NN, \ \ \ \ n\to\IF,
\EQNY
or equivalently, for any $x\in\R,$
\BQNY
\lim_{n\to\IF}\pk{N_n(x)<k} = \int_{-\IF}^\IF \exp\LT(-e^{-x + y \sigma_0c\beta/H }\RT) \sum_{l=0}^{k-1}\frac{e^{l (-x+y \sigma_0c\beta/H)}}{l !}  \varphi(y) dy,
\EQNY
where $\varphi(y)=(2\pi)^{-1/2} e^{-y^2/2}, y\in\R,$ is the density function of the standard Normal distribution. That is, the number of exceedances $N_n(x)$ is approximately mixed Poisson distributed with random intensity $\lambda(x) = e^{-x + \NN \sigma_0 c\beta/H}$.
\end{itemize}

\ET

 \begin{remarks} \label{Rem:1}
(a). In the case $\beta>2H-H_0$, it can be understood as  that the dependence among $\{Q_i\}_{i\ge1}$ is so strong  that in the limit the sequence will have either infinitely many or no exceedances of a high-level. 

(b). It is interesting to notice that the above three types of limiting result for $k=1$  (i.e., Normal, Gumbel and a mixture of them) 
resemble the classical results for the stationary Normal sequences. 

(c). It is worth noting that in the case of stationary Normal sequence (and many other general stationary sequences) it is the  the Poisson (or Cox) point process convergence that is first obtained which implies the convergence for order statistics. As discussed in the Introduction this approach might not work here due to  non-existence of a mixing condition for the stationary sequence $\{Q_i\}_{i\ge1}$. 
Here we directly  prove a weak convergence result for the order statistics which is equivalent to a (mixed) Poisson distribution convergence under the last two scenarios.

(d). A stationary sequence $\{\xi_i\}_{i\ge 1}$ is said to admit a phantom distribution function $G$ if
$$
\pk{\max_{i\le n} \xi_i \le  u_n  }- G^n(u_n) \to 0, \ \ \ \ n\to\IF,
$$ 
for every sequence $\{u_n\}_{n\ge 1}\subset \R,$ see e.g., \cite{DJL15} and references therein.
It can be shown that under scenarios (i) and (iii) the stationary sequence $\{Q_i\}_{i\ge1}$ does not admit a phantom distribution function, whereas under scenario (ii) it admits a continuous phantom distribution function. 
A  proof is given in Section 4 following the proof of Theorem \ref{Thm:1}.

(e).  It is of interest to study the existence and value (if exists) of extremal index $\theta$ of the  stationary sequence $\{Q_i\}_{i\ge1}$; see Section 3.7 of \cite{leadbetter1983extremes} for a definition of extremal index. To this end, the asymptotics of $\pk{Q_1> u}$, as $u\to\IF$, seems to be a key tool; some results regarding this asymptotics have been obtained in 
\cite{HS06} under some additional conditions (see A1 and A2 therein for slightly general Gaussian processes) which are assumed to hold here for simplicity. In order to save some space we only give some comments, omitting technical assumptions and derivations, for this remark. We can show that under scenario (ii), the extremal index $\theta =1$. This can be checked by choosing $u_n(x)=a_n x+b_n$ and using the asymptotics of Theorem 2.1 combined with formulas (5) and (7) in \cite{HS06}. In fact, it is quite intuitive that when $H$ is large enough (in the sense of scenario (ii)) the stationary sequence $\{Q_i\}_{i\ge1}$ shows a strong independence which allows it to have an associated independent sequence in the sense of \cite{leadbetter1983extremes} and thus $\theta=1$. Similarly, we can check that the extremal index does not seem to make sense under scenario (iii), this is understandable due to the mixture type of the limiting distribution in (iii) of Theorem \ref{Thm:1}. Finally, under scenario (i) we conjecture that $\theta=0$, this is understandable intuitively due to some strong clustering property discussed in remark (a) above. It seems hard to confirm such a result in general because of the complicated higher than first order asymptotics for the function $f_u^2(s)$, as $u\to\IF,$ in (5)  of \cite{HS06} under this scenario. However, we can easily verify this conjecture for the Brownian model with linear drift, using explicit formulas. 


(f). We remark that extensions of Theorem \ref{Thm:1} to multivariate order statistics of the form 
\BQNY
\vk M_n^{(k)} :=\LT(\max_{i\le n}^{(k)}   Q_{1,i}, \ \max_{i\le n}^{(k)}   Q_{2,i}, \cdots, \ \max_{i\le n}^{(k)}   Q_{d,i}\RT),
\EQNY 
can be done similarly, where $Q_{l,i} = \sup_{t\ge 0}( X_i^{(l)}(t)+\sigma_0 X^{(l)}(t)-ct^\beta)$ with $\{X_i^{(l)}(t), t\ge0\}, {l=1,\cdots d, i=1,\cdots, n}$ being independent copies of a self-similar Gaussian process and $\{(X^{(1)}(t), \cdots, X^{(d)}(t)), t\ge0\}$ being a $d$-dimensional  self-similar Gaussian process. We refer to 
\cite{ACLP13} for examples of multivariate self-similar Gaussian processes which include some multivariate fBm as special case.

\end{remarks}

\COM{
\begin{remark} We remark on a possible extension of Theorem \ref{Thm:1} to multivariate case. To ease the notation and to illustrate the idea, we shall consider a particular case.
 Let $\{X_i(t), t\ge0 \}, i =1,2,\ldots,$ and $\{Y_i(t), t\ge0 \}, i =1,2,\ldots,$ be two independent sequences of self-similar Gaussian processes with a common self-similarity index $H\in(0,1)$. Further, let $(X(t), Y(t)), t\ge0$ be an independent  two-dimensional self-similar Gaussian process  with self-similarity index $H_0\in(0,1)$, i.e., for any $a>0$,
 \BQNY
\{ (X(at), Y(at))\}_{t\ge0} \overset{d} = a^{H_0} \{ (X(t), Y(t))\}_{t\ge0}.
 \EQNY
This is a particular definition of the known  operator-self-similar processes (see, e.g., \cite{ACLP13}) and includes some multivariate fBm as special case. Define
\BQNY
\vk M_n^{(k)} :=\LT(\max_{i\le n}^{(k)}   Q_{1,i}, \ \max_{i\le n}^{(k)}   Q_{2,i}\RT),
\EQNY
with $Q_{1,i} = \sup_{t\ge 0}( X_i(t)+X(t)-ct^\beta)$ and $Q_{2,i} = \sup_{t\ge 0}( Y_i(t)+Y(t)-ct^\beta)$ for some $c>0$ and $\beta>\max(H,H_0)$. Using similar arguments as in the proof of Theorem \ref{Thm:1}, it can be shown that
\begin{itemize}
\item[(i).] If $\beta>2H-H_0$, then

\BQNY 
t_0^{-H_0}  b_n^{- H_0 /\beta}  \LT(\vk M_n^{(k)} -  (b_n, b_n) \RT)   \ \overset{d}\to \  (\NN_1, \NN_2), \ \ \ \ n\to\IF.
\EQNY

\item[(ii).] If   $\beta<2H-H_0$, then

\BQNY
a_n^{-1}\LT(\vk M_n^{(k)} -  (b_n, b_n) \RT)   \ \overset{d}\to \  (\Lambda_1^{(k)},\Lambda_2^{(k)}), \ \ \ \ n\to\IF.
\EQNY

\item[(iii).] If   $\beta=2H-H_0$, then

\BQNY
a_n^{-1}\LT(\vk M_n^{(k)} -  (b_n, b_n) \RT) \ \overset{d}\to \  \LT(\Lambda_1^{(k)} + \frac{c\beta}{H} \NN_1, \Lambda_2^{(k)} + \frac{c\beta}{H} \NN_2\RT), \ \ \ \ n\to\IF.
\EQNY

\end{itemize}
In the above, $\Lambda_1^{(k)}, \Lambda_2^{(k)}$ are two independent copies of $\Lambda^{(k)}$, $(\NN_1, \NN_2)$ is an independent of $\Lambda_1^{(k)}$ and  $\Lambda_2^{(k)}$   Normal random vector which is equal in distribution to $(X(1), Y(1))$.

\end{remark}
}

The next result shows that for the stationary sequence $\{Q_i\}_{i\ge 1}$, the (absolute) moments of the normalised order statistics
converge to the  (absolute) moments of the corresponding limiting   random variable. 


\BT \label{Thm:2}

Under the assumptions of Theorem \ref{Thm:1}, we have, for any $\lambda>0,$
\begin{itemize}
\item[(i).] If $\beta>2H-H_0$, then

\BQNY 
\lim_{n\to\IF}\E{\LT|  \sigma
_0^{-1}t_0^{-H_0}  b_n^{ -H_0 /\beta}\LT( M_n^{(k)} -  b_n \RT) \RT|^\lambda  } = \E{\abs{\NN}^\lambda}.
\EQNY

\item[(ii).] If  $\beta<2H-H_0$, then
\BQNY
\lim_{n\to\IF}\E{\LT| a_n^{-1}\LT(M_n^{(k)} -  b_n\RT)  \RT|^\lambda  } = \E{ \abs{ \Lambda^{(k)}}^\lambda}.
\EQNY

\item[(iii).] If  $\beta=2H-H_0$, then

\BQNY
\lim_{n\to\IF}\E{\LT|  a_n^{-1}\LT(M_n^{(k)} -  b_n\RT) \RT|^\lambda  } = \E{ \abs{  \Lambda^{(k)} + \frac{\sigma_0c\beta}{H} \NN}^\lambda}.
\EQNY
\end{itemize}

\ET
\begin{remark}
We can see  from the proof of Theorem \ref{Thm:2} that, when $\lambda$ is an integer, the above convergence results still hold for moments without the modulus. Absolute moments of the limiting distributions can sometimes be given more explicitely, for example, it follows from   \cite{Win14}  that
 $\E{\abs{\NN}^\lambda}= \frac{2^{\lambda/2}}{\sqrt \pi} \Gamma\LT(\frac{\lambda+1}{2}\RT),$
 with $\Gamma(\cdot)$ the Gamma function. Furthermore, by a change a variable formula, we can obtain $\E{ \abs{ \Lambda^{(1)}}^\lambda}=\int_0^\IF \abs{\log y}^\lambda e^{-y} dy.$ The formula for other distributions seems to be complicate and thus omitted here. Moreover, these moments can be easily approximated by using Monte Carlo simulations.
\end{remark}
\COM{
\begin{remark} Denote
$\wM_n^{(k)}= \max^{(k)}_{i\le n} \wQ_i$ with the IID sequence $\wQ_i, i=1,2,\ldots$ defined in \eqref{eq:wMQ}. As a by-product, we can show that
\BQN\label{eq:wMab}
\lim_{n\to\IF}\E{\LT| a_n^{-1} \LT(\wM_n^{(k)} -  b_n\RT)\RT|^\lambda  } = \E{ \abs{ \Lambda^{(k)}}^\lambda}.
\EQN
As far as we are aware, general results for $k$th largest order-statistics of IID Weibull-like random variables of this type have been discussed in the literature. This is a problem that is of independent interest. Similar results exist only when $k=1$, that is, for the maximum (see, e.g., \cite{Pick68} or \cite{Res87}).
We include a short of proof for \eqref{eq:wMab} in the Appendix.
\end{remark}
}

The following theorem presents analogues of Theorem \ref{Thm:1} and Theorem \ref{Thm:2} for a non-stationary sequence $\{Q_i\}_{i\ge 1}$ with general $c_i$'s.

\BT \label{Thm:1-1}
Let $M_n^{(k)}, n\in \mathbb N$ be defined in  \eqref{eq:Qi}-\eqref{eq:MQ} with $\sigma=1$,  and \eqref{eq:cc} holds with some $m_n < n$ such that $\lim_{n\to\IF}m_n/n = p\in(0,1]$,  and let $b_n, a_n, n\in \mathbb N$ be given as in
\eqref{eq:bnan}.
 Assume $\beta>\max(H,H_0)$. Then, the claims of (i)-(iii) in Theorem \ref{Thm:1} and Theorem \ref{Thm:2} hold true when replacing  $a_n, b_n$ with $a_{m_n}, b_{m_n}$, respectively.
\ET

\begin{remark}
The above result is understandable 
intuitively as follows: The probability of exceeding a high-level threshold by $Q_i,$ for any $i> m_n$ is much less than that of $Q_i, i\le m_n$,  so a lower threshold $u_{m_n}(x)$ (defined through $m_n$ instead of $n$) is needed in order to derive the same limiting distribution as for the stationary case. In this sense, the above results for the non-stationary sequence $\{Q_i\}_{i\ge 1}$ can be seen as a thinning version of  the results  in Theorem \ref{Thm:1} and Theorem \ref{Thm:2}. 
\end{remark}



We conclude this section with an example, where we derive corresponding results for the  fBm  model  with a linear drift (i.e., $\beta=1$). For a sfBm $\{X_1(t), t\ge 0\}$ with  Hurst index $H\in (0,1)$, 
\BQNY
\Cov(X_1(t),X_1(s))=\frac{1}{2}(t^{2H}+s^{2H}-\mid t-s\mid^{2H}),\ \ \ t,s\ge0.
\EQNY
\Pe{One can check that sfBm $X_1$ fulfills \eqref{eq:locstaXH}} with $K(t)=t_0^{-H} t^H = (H/(c(1- H)) )^{- H/\beta} t^H, t\ge0$. 
\Pe{Thus, by Proposition \ref{prop1}, we have
\BQNY
\pk{\sup_{t\geq0} (X_1(t)-ct) >u} = R(u)\exp\LT(-\frac{u^{\tau}}{2A^2}\RT)(1+o(1))
\EQNY
as $u\to\IF$, where}
\BQN\label{eq:ATR}
A=\frac{H^H (1-H)^{1-H}}{c^H}, \ \ \tau =2(1-H), \ \ R(u)=  \frac{ 2^{-\frac{1}{2H}}  \mathcal{H}_{2H}}{\sqrt{H(1-H)}}\left(\frac{c^H u^{1-H}}{H^H (1-H)^{1-H}}\right)^{\frac{1}{H}-2}.
\EQN


\begin{korr} \label{Cor:fBm}
Let $\{X_i(t), t\ge 0\}, i=1,2,\ldots,$  be independent sfBm's with common Hurst index $H\in (0,1)$ and $\{X(t), t\ge 0\}$  be another independent sfBm with Hurst index $H_0\in (0,1)$. Assume $\sigma=1$ and  $\beta=1$. Then, the claims in Theorems \ref{Thm:1}, \ref{Thm:2} and  \ref{Thm:1-1} are valid, with $b_n, a_n$ in \eqref{eq:bnan} defined through \eqref{eq:ATR}.
\end{korr}
\begin{remark}
Particularly, if $H=H_0=1/2$, $c=c_i,i\ge1$ and $k=1$, we recover Theorem 5.2 of \cite{MSJVZ21}.
\end{remark}

\section{ Further results and Proofs} \label{Sec:proof}

Before starting the proof, we first give some auxiliary results which will be used later. Recall the key point $t_0=\LT(\frac{H}{c(\beta-H)}\RT)^{\frac{1}{\beta}}$ as given above the formula \eqref{eq:AB}.
The  lemma below is about limiting properties of $a_n$ and $b_n$ which can be obtained immediately from their definition.

\BEL
For $a_n$ and $b_n$ in \eqref{eq:bnan}, we have
\BQN  \label{limit-abn}
\lim_{n\to\IF }\frac{b_n^{H_0/\beta}}{a_n}  = \tau (2A^2)^{\frac{H_0-\beta}{2(\beta-H)}} \lim_{n\to\IF} (\log n)^{\frac{\beta-2H+H_0}{2(\beta-H)}}
 =  \left\{
\begin{array}{ll}
\IF, &  \mbox{if $\beta> 2H-H_0$,}\\[0.1cm]
0, & \mbox{if $\beta< 2H-H_0$,}\\[0.1cm]
\tau /(2A^2),  & \mbox{if $\beta= 2H-H_0$.}
\end{array}
\right.
\EQN
Furthermore, 
as $n\to\IF,$
\BQN \label{tayex-bn}
b_n^{1-H/\beta}
 = A\sqrt{2\log n} \LT( 1+ \frac{1}{2} (\log n)^{-1} \log(R((2A^2\log n)^{\frac{\beta}{2(\beta-H)} }))\oo      \RT) . 
\EQN
\EEL


\BEL\label{Lem:hatc}
For any $\vp_0\in(0,t_0)$, there exists some small \Pe{$\hat c\in(0,\min(1,c))$} such that
\BQN\label{def-sigsta}
\sigma_{\ast}=\sigma_{\ast}(\vn_0, \hat c)  :=  \max\LT\{\max_{t\in[0,t_0-\vp_0]} \frac{t^H}{1-{\hat c} +ct^\beta},\
 \max_{t\geq t_0+\vp_0} \frac{t^H}{1-{\hat c} +(c-\hat c)t^\beta} \RT\}
 <  \frac{t_0^H}{1+ct_0^\beta}=A.
\EQN

\EEL
\prooflem{Lem:hatc} We only show the proof for the maximum taken over $[0,t_0-\vp_0]$, since similar arguments  also apply   to the  second maximum taken over $[t_0+\vp_0, \IF)$.
Note that, for any $\vp_1>-1$ and $\vp_2\in[0,c)$,
\BQN\label{maxpointre}
\mathop{\arg\max}\limits_{t\geq0}\frac{t^{H}}{1+\vp_1+(c-\vn_2)t^{\beta}}=t_0\LT( \frac{c}{c-\vn_2} (1+\vp_1)\RT)^{1/\beta}.
\EQN
[This formula is given in a general form which is also helpful for later.]
 Thus, for any $\vp_0\in(0,t_0)$, we can find some small enough $\hat c$ such that
 \BQNY
\max_{t\in[0,t_0-\vp_0]} \frac{t^H}{1-{\hat c} +ct^\beta}
 =   \frac{(t_0-\vp_0)^H}{1-{\hat c} +c(t_0-\vp_0)^\beta}< \frac{t_0^H}{1+ct_0^\beta},
\EQNY
 where the last inequality follows by  \eqref{maxpointre} with $\vp_1=\vp_2=0$. 
 This completes the proof. \QED

\subsection{Proof of Theorem \ref{Thm:1}}  In the following subsections, we first present the proof for scenario (i) and then a generic proof for scenarios (ii)-(iii).
\subsubsection{Proof for (i) }
We need to show that, for any $x\in \R$,
\BQNY
\pk{ b_n^{-H_0/\beta} \LT(M_n^{(k)} -b_n\RT) >x} \to \pk{\sigma_0 t_0^{H_0} \NN >x},\ \ \ \ n\to\IF.
\EQNY
\Ji{We will consider asymptotic lower and upper bounds, respectively.}
First, we have, from Lemma \ref{Lem:1} below, that
\BQN
\pk{b_n^{-H_0/\beta}  \LT(M_n^{(k)} -b_n\RT) >x} &\ge& \pk{ b_n^{-H_0/\beta} \LT(\overset{(k)} {\underset{i\le n}{\max} }\ X_i (t_0 b_n^{1/\beta}) +\sigma_0 X(t_0 b_n^{1/\beta}) -c (t_0 b_n^{1/\beta})^\beta -b_n \RT)  >x}\nonumber\\
&\to&\pk{\sigma_0t_0^{H_0} \NN >x},\ \ \ \ n\to\IF, \label{eq:t0NN}
\EQN
which yields the required lower bound. Next,  for any $\vn_0\in(0,t_0)$, we  introduce the following notation,
\BQNY
&& A_{1,i}=\LT\{ b_n^{-H_0/\beta}  \LT(\sup_{(t_0-\vn_0) b_n^{1/\beta}\le t\le (t_0+\vn_0) b_n^{1/\beta}} (X_i (t) +\sigma_0X(t)-c t^\beta) -b_n \RT)    >x\RT\}, \ \ 1\le i\le n,\\
&&A_{2,i}=\LT\{ b_n^{-H_0/\beta} \LT(\sup_{0\le t\le (t_0-\vn_0) b_n^{1/\beta}} (X_i (t) +\sigma_0X(t)-c t^\beta) -b_n \RT)  >x\RT\},\ \ 1\le i\le n,\\
&&A_{3,i}=\LT\{ b_n^{-H_0/\beta} \LT(\sup_{t\ge (t_0+\vn_0) b_n^{1/\beta}}(X_i (t)+\sigma_0X(t) -c t^\beta) -b_n \RT)  >x\RT\},\ \ 1\le i\le n,\\
&& B_2= \cup_{i\le n} A_{2,i}, \ \ \  \ B_3= \cup_{i\le n} A_{3,i}.
\EQNY
We derive that
\BQN\label{eq:ABB}
\pk{b_n^{-H_0/\beta}  \LT(M_n^{(k)} -b_n\RT) >x} &=&\pk{  \underset{ \subseteq \{1,\ldots,n\} } {\underset{\{ j_1,\ldots, j_k\}} \bigcup }    \Big(  \LT(\cup_{i=1}^3 A_{i,j_1}\RT) \cap \cdots \cap \LT(\cup_{i=1}^3 A_{i,j_k}\RT)     \Big)  }\nonumber\\
&=& \pk{  \underset{ \subseteq \{1,\ldots,n\} } {\underset{\{ j_1,\ldots, j_k\}} \bigcup }    \Big(  \LT(\cup_{i=1}^3 A_{i,j_1}\RT) \cap \cdots \cap \LT(\cup_{i=1}^3 A_{i,j_k}\RT) \Big), (B_2 \cup B_3)^c    }\nonumber\\
&&+\, \pk{  \underset{ \subseteq \{1,\ldots,n\} } {\underset{\{ j_1,\ldots, j_k\}} \bigcup }    \Big(  \LT(\cup_{i=1}^3 A_{i,j_1}\RT) \cap \cdots \cap \LT(\cup_{i=1}^3 A_{i,j_k}\RT) \Big), (B_2 \cup B_3)      }\\
& \le& \pk{\underset{ \subseteq \{1,\ldots,n\} } {\underset{\{ j_1,\ldots, j_k\}} \bigcup }    \Big(   A_{1,j_1} \cap \cdots \cap   A_{1,j_k}  \Big)} +  \pk{B_2 } + \pk{B_3},\nonumber
\EQN
where ${\underset{\{ j_1,\ldots, j_k\} \subseteq\{1,\ldots,n\} } \bigcup } $ denotes the union of all the possible combinations of $j_1,\ldots, j_k$ drawn without replacement from $\{1,\ldots,n\}$. Thus, the above inequality can be re-written as
\BQNY 
\pk{b_n^{-H_0/\beta}  \LT(M_n^{(k)} -b_n\RT) >x} &\le& \pk{ \max_{i\le n}^{  (k)} \sup_{(t_0-\vn_0) b_n^{1/\beta}\le t\le (t_0+\vn_0) b_n^{1/\beta}}\frac{X_i (t) +\sigma_0X(t) -c t^\beta -b_n}{ b_n  ^{H_0/\beta}} >x }\nonumber\\
&&  +\ \pk{ \max_{i\le n} \sup_{0\le t\le (t_0-\vn_0) b_n^{1/\beta}}\frac{X_i (t) +\sigma_0X(t) -c t^\beta -b_n}{ b_n  ^{H_0/\beta}} >x } \\
&&+\ \pk{ \max_{i\le n} \sup_{ t\ge (t_0+\vn_0) b_n^{1/\beta}}\frac{X_i (t) +\sigma_0X(t) -c t^\beta -b_n}{ b_n  ^{H_0/\beta}} >x }. \nonumber
\EQNY
In view of \eqref{negpart1} and \eqref{negpart2} in Lemma \ref{Lem:2} below, we know that the last two terms on the right-hand side converge to 0, as $n\to\IF$. For the remaining first term, it follows, by self-similarity, that
\BQNY
&&\quad\pk{ \overset{(k)} {\underset{i\le n}{\max} }\sup_{(t_0-\vn_0) b_n^{1/\beta}\le t\le (t_0+\vn_0) b_n^{1/\beta}}\frac{X_i (t) +\sigma_0X(t) -c t^\beta -b_n}{ b_n  ^{H_0/\beta}} >x }\nonumber\\
&&  \le \pk{ \overset{(k)} {\underset{i\le n}{\max} }\ \sup_{t\ge 0}\frac{X_i (t)   -c t^\beta -b_n}{ b_n  ^{H_0/\beta}} + \sup_{(t_0-\vn_0) b_n^{1/\beta}\le t\le (t_0+\vn_0) b_n^{1/\beta}}\frac{\sigma_0X(t)}{b_n  ^{H_0/\beta}}>x } \\
&&  = \pk{   b_n  ^{-H_0/\beta} (\wM_n^{(k)} -b_n ) + \sup_{1\le t\le (t_0+\vn_0) /(t_0-\vn_0)} \sigma_0X(t) (t_0-\vn) ^{H_0} >x}\nonumber
\EQNY
with $\wM_n^{(k)}$ defined in \eqref{def:MteQte} with $c=c_i, i\ge1$ in the stationary case. 
Therefore, we derive from \eqref{eq:wMb} and \eqref{limit-abn} that, for $\beta>2H-H_0$,
\BQNY
 b_n  ^{-H_0/\beta} (\wM_n^{(k)} -b_n )  = \frac{a_n}{b_n  ^{H_0/\beta}}\frac{\wM_n^{(k)} -b_n}{ a_n } \ \overset{d}\to \ 0, \quad n\to\IF,
\EQNY
and thus
\BQNY
&&\quad \Pe{\lim_{\vp_0\to0}}\limsup_{n\to\IF} \pk{ \overset{(k)} {\underset{i\le n}{\max} }\sup_{(t_0-\vn_0) b_n^{1/\beta}\le t\le (t_0+\vn_0) b_n^{1/\beta}}\frac{X_i (t) +\sigma_0X(t) -c t^\beta -b_n}{ b_n  ^{H_0/\beta}} >x }\\
&&\leq \lim_{\vn_0\to 0 }\pk{  \sup_{1\le t\le (t_0+\vn_0) /(t_0-\vn_0)} \sigma_0 X(t) (t_0-\vn_0) ^{H_0}>x} = \pk{\sigma_0t_0^{H_0} \NN >x},
\EQNY
%
which gives the required upper bound. The proof is complete.  \QED

Below we present the two lemmas used in this proof.

\BEL \label{Lem:1}
Under the assumption of Theorem \ref{Thm:1} \Ji{and the condition in (i) (i.e., $\beta>2H-H_0$),} we have, \Ji{as $n\to\IF,$}
\BQNY
b_n  ^{-H_0/\beta}\LT( \overset{(k)} {\underset{i\le n}{\max} }\ X_i (t_0 b_n^{1/\beta})) +\sigma_0X(t_0 b_n^{1/\beta}) -c (t_0 b_n^{1/\beta})^\beta -b_n \RT)  \overset{d} \to \  \sigma_0t_0^{H_0} X(1). 
\EQNY

\EEL

\prooflem{Lem:1} First, by self-similarity, 
\BQNY
&&\quad b_n  ^{-H_0/\beta}\LT( \overset{(k)} {\underset{i\le n}{\max} }\ X_i (t_0 b_n^{1/\beta})) +\sigma_0X(t_0 b_n^{1/\beta}) -c (t_0 b_n^{1/\beta})^\beta -b_n \RT) \\
&&\overset{d}  = \ t_0^H \LT( \overset{(k)} {\underset{i\le n}{\max}} X_i (1 )   -\frac{1+c t_0 ^\beta}{ t_0^H}  b_n^{1-H/\beta} \RT) b_n  ^{(H-H_0)/\beta}
 + \sigma_0 t_0^{H_0} X(1). 
\EQNY
 It is sufficient to show that
\BQN\label{eq:Xito0}
\LT(\overset{(k)} {\underset{i\le n}{\max} }\ X_i (1 )   -\frac{1+c t_0 ^\beta}{ t_0^H}  b_n^{1-H/\beta} \RT) b_n  ^{(H-H_0)/\beta} \overset{d}\to 0, \ \ \ \ n\to\IF.
\EQN
For the IID standard Normal sequence $X_i(1), i=1,2,\ldots$, 
we have from Proposition \ref{Prop:GamLim} (see also Theorem 1.5.3 in \cite{leadbetter1983extremes}) that
\BQN\label{eq:NLam}
\Ji{Z_n^{(k)}:=}\sqrt{2 \log n}\LT(\overset{(k)} {\underset{i\le n}{\max} }\ X_i(1) - \LT(\sqrt{2 \log n} - \frac{\log(4\pi \log n)}{2 \sqrt{2 \log n}}\RT)\RT) \overset{d} \to \Lambda^{(k)}, \ \ \ \ n\to\IF.
\EQN
Then, we can rewrite the left-hand side of \eqref{eq:Xito0} as
\BQN \nonumber
&&\quad\LT(\overset{(k)} {\underset{i\le n}{\max} }\ X_i (1 )   -\frac{1+c t_0 ^\beta}{ t_0^H}  b_n^{1-H/\beta} \RT) b_n  ^{(H-H_0)/\beta}\\ \nonumber
&& = \frac{\sqrt{2 \log n} \LT(\overset{(k)} {\underset{i\le n}{\max} }\ X_i (1 ) -\LT(\sqrt{2 \log n} - \frac{\log(4\pi \log n)}{2 \sqrt{2 \log n}}\RT)  \RT)}{\sqrt{2 \log n}\  b_n  ^{(H_0-H)/\beta}}
- \frac{ \frac{1+c t_0 ^\beta}{ t_0^H}  b_n^{1-H/\beta} - \LT(\sqrt{2 \log n} - \frac{\log(4\pi \log n)}{2 \sqrt{2 \log n}}\RT)  }{b_n  ^{(H_0-H)/\beta}}\\ \label{nor-asy}
&& =:  \frac{Z_n^{(k)}}{\sqrt{2 \log n}\  b_n  ^{(H_0-H)/\beta}} -r_n.
\EQN
By the definition of $b_n$ in \eqref{eq:bnan} \Pe{and the assumption $\beta>2H-H_0$}, we have
\Pe{
\BQN\label{lim-logn-bn}
\quad\lim_{n\to\IF}\sqrt{2\log n}\  b_n  ^{(H_0-H)/\beta}
=\lim_{n\to\IF}A^{\frac{H_0-H}{\beta-H}}  (2\log n ) ^{\frac{\beta+H_0-2H}{2(\beta-H)}} \LT(1+ o(1)
 \RT)^{(H_0-H)/\beta} =\IF,
\EQN
}
\Pe{which together with} Taylor expression \eqref{tayex-bn} implies that
\Pe{
\BQN\label{asy-bnH-H0}
\lim_{n\to\IF}r_n=\frac{1}{\sqrt{2\log n}\ b_n  ^{(H_0-H)/\beta}} \LT(   \log(R((2A^2\log n)^{1/\tau}))(1+o(1)) + 2^{-1} \log(4\pi \log n)  \RT)=0.
\EQN
}
Consequently, substituting \eqref{eq:NLam} and \eqref{lim-logn-bn}-\eqref{asy-bnH-H0} into \eqref{nor-asy},
we get \eqref{eq:Xito0}.
This completes the proof.  \QED

\BEL \label{Lem:2} Under the assumption of Theorem \ref{Thm:1}, we have, \Ji{for any $\vn_0\in(0,t_0)$ and any $x\in\R$,}
\BQN\label{negpart1}
\quad \lim_{n\to\IF}\pk{ \max_{i\le n} \sup_{0\le t\le (t_0-\vn_0) b_n^{1/\beta}}\frac{X_i (t) +\sigma_0X(t) -c t^\beta -b_n}{ b_n  ^{H_0/\beta}} >x }=0,
\EQN
\BQN\label{negpart2}
\lim_{n\to\IF}\pk{ \max_{i\le n} \sup_{t\ge (t_0+\vn_0) b_n^{1/\beta}}\frac{X_i (t) +\sigma_0X(t) -c t^\beta -b_n}{ b_n  ^{H_0/\beta}} >x }=0.
\EQN
\EEL

\prooflem{Lem:2} We first prove \eqref{negpart1}. Note, by self-similarity,
\BQN
&&\quad\pk{ \max_{i\le n} \sup_{0\le t\le (t_0-\vn_0) b_n^{1/\beta}}\frac{X_i (t) +\sigma_0X(t) -c t^\beta -b_n}{ b_n  ^{H_0/\beta}} >x }\nonumber \\
&&\le   \pk{ \max_{i\le n} \sup_{0\le t\le (t_0-\vn_0) b_n^{1/\beta}}\frac{X_i (t)  - c t^\beta -b_n}{ b_n  ^{H
_0/\beta}} +\sup_{0\le t\le (t_0-\vn_0)b_n  ^{1/\beta} } \frac{\sigma_0X(t)}{b_n  ^{H_0/\beta}} > x }  \label{eq:UUU} \\
&&=   \pk{ \max_{i\le n}\sup_{0\le t\le (t_0-\vn_0) } \LT( b_n  ^{(H-H_0)/\beta} X_i (t)  - (1+c t^\beta )b_n  ^{1-H_0/\beta}  \RT)+\sup_{0\le t\le (t_0-\vn_0) }\sigma_0X(t) > x }. \nonumber
\EQN
Since $\sup_{0\le t\le (t_0-\vn_0) }X(t)<\IF$ \Ji{a.s.}, it is sufficient to show that, for any  $x\in\R$,
\BQN\label{eq:Jnx10}
J_1(n,x):=\pk{ \max_{i\le n}\sup_{0\le t\le (t_0-\vn_0)}\LT( b_n  ^{(H-H_0)/\beta} X_i (t)  - (1+c t^\beta )b_n  ^{1-H_0/\beta} \RT) > x } \to 0,
\EQN
as $n\to\IF$. \Pe{For the fixed $\vn_0$, choosing a small $\hat c\in(0,1)$ satisfying \eqref{def-sigsta}, then using  Borell-TIS inequality (see, e.g., Theorem 2.1.1 in \cite{AdlerTaylor}), we have, for large enough $n$ such that
$xb_n^{H_0/\beta-1}> - \hat c$,}
\BQNY
J_1(n,x)
&\le & n \ \pk{  \sup_{0\le t\le (t_0-\vn_0) } \frac{X_i (t)}{  1+c t^\beta   + x   b_n  ^{H_0/\beta-1}} >b_n  ^{1-H/\beta}  }\\
&\le & n \ \pk{  \sup_{0\le t\le (t_0-\vn_0) } \frac{X_i (t)}{  1 -\hat c +c t^\beta } >b_n  ^{1-H/\beta}  }\\
&\le & n \exp\LT( -( b_n  ^{1-H/\beta} - K_1)^2 \LT( \sup_{t\in[0,t_0-\vp_0]} \frac{2t^{2H}}{  \LT(1
-\hat c + c  t^\beta\RT)^2}\RT)^{-1} \RT)\\
 &\le&  \exp\LT( -\frac{( b_n  ^{1-H/\beta} - K_1)^2}{2 \sigma_*^2} + \log n\RT),
\EQNY
where $K_1:=  \E{\sup_{0\le t\le (t_0-\vn_0) } X_1 (t)/(1 -\hat c +c t^\beta) }< b_n  ^{1-H/\beta}$ \Ji{for all large enough $n$}, and the last inequality follows from \eqref{def-sigsta}. Furthermore, by \eqref{tayex-bn} \Ji{and  \eqref{def-sigsta}} we have
\BQNY
\lim_{n\to\IF} \frac{ b_n  ^{2-2H/\beta} }{2 \sigma_*^2 \log n} >1,
\EQNY
implying \eqref{eq:Jnx10}.
Thus, \eqref{negpart1} is established.  Next, by a similar argument we derive, 
for some $\hat c\in(0,\min(1,c))$ satisfying \eqref{def-sigsta}, that
\BQNY
&&\quad\pk{ \max_{i\le n} \sup_{t\ge (t_0+\vn_0) b_n^{1/\beta}}\frac{X_i (t) +\sigma_0X(t) -c t^\beta -b_n}{ b_n  ^{H_0/\beta}} >x }\\
&& \le\pk{ \max_{i\le n} \sup_{ t\ge (t_0+\vn_0) b_n^{1/\beta}}\frac{X_i (t)  - (c-\hat c) t^\beta -b_n}{ b_n  ^{H_0/\beta}} +  \sup_{ t\ge (t_0+\vn_0) b_n^{1/\beta}}\frac{\sigma_0X(t)  - \hat c  t^\beta }{ b_n  ^{H_0/\beta}} >x } \\
&&\le  \pk{ \max_{i\le n} \sup_{ t\ge (t_0+\vn_0) }\LT(b_n  ^{(H-H_0)/\beta}X_i (t)  - (1+(c-\hat c) t^\beta) b_n  ^{1-H_0/\beta}\RT) + b_n  ^{-H_0/\beta}\sup_{t\geq0}\LT( \sigma_0X(t)-\hat c t^\beta\RT) >x}.
\EQNY
Note that 
$\lim_{n\to\IF}b_n  ^{-H_0/\beta}\sup_{t\geq0} (\sigma_0X(t)-\hat c t^\beta)  =0$ a.s.. Thus, in order  to prove \eqref{negpart2}, it is sufficient to show that, for any fixed $x\in\R$,
\BQNY
J_2(n,x):=\pk{ \max_{i\le n} \sup_{ t\ge (t_0+\vn_0) }\LT(b_n  ^{(H-H_0)/\beta}X_i (t)  - (1+(c-\hat c) t^\beta) b_n  ^{1-H_0/\beta}\RT)  >x} \to 0,
\EQNY
as $n\to\IF$. Again, by Borell-TIS inequality, we have, 
\BQNY
J_2(n,x)
&\le & n \ \pk{  \sup_{t\geq (t_0+\vn_0)} \frac{X_i (t)}{  1+(c-\hat c) t^\beta   + x   b_n  ^{H_0/\beta-1}} >b_n  ^{1-H/\beta}  }\\
&\le & n \ \pk{  \sup_{t\geq (t_0+\vn_0)} \frac{X_i (t)}{  1 -\hat c +(c-\hat c) t^\beta   } >b_n  ^{1-H/\beta}  }\\
&\le & n \exp\LT( -( b_n  ^{1-H/\beta} - K_2)^2 \LT( \sup_{t\geq (t_0+\vn_0)} \frac{2t^{2H}}{  \LT(1
-\hat c + (c-\hat c)  t^\beta\RT)^2}\RT)^{-1} \RT)\\
 &\le&  \exp\LT( -\frac{( b_n  ^{1-H/\beta} - K_2)^2}{2 \sigma_*^2} + \log n\RT)\to 0 , \quad n\to\IF,
\EQNY
where $K_2:=  \E{\sup_{ t\ge (t_0+\vn_0) } \frac{X_1 (t) }{ 1+(c-\hat c ) t^\beta}}< \IF$.
Thus,  the proof is complete. \QED

\medskip

\Pe{Before proving scenarios (ii) and (iii), we shall  derive two important lemmas below.}

\BEL \label{Lem:negli} Under the assumptions of Theorem \ref{Thm:1} \Ji{and the conditions in (ii) and (iii)} (i.e., $\beta\leq 2H-H_0$), we have, for any $\vn_0\in(0,t_0)$ and any $x\in\R$,
\BQN\label{negpart-ii-1}
\quad \lim_{n\to\IF}\pk{ \max_{i\le n} \sup_{0\le t\le (t_0-\vn_0) b_n^{1/\beta}}\frac{X_i (t) +\sigma_0X(t) -c t^\beta -b_n}{ a_n} >x }=0,
\EQN
\BQN\label{negpart-ii-2}
\lim_{n\to\IF}\pk{ \max_{i\le n} \sup_{t\ge (t_0+\vn_0) b_n^{1/\beta}}\frac{X_i (t) +\sigma_0X(t) -c t^\beta -b_n}{ a_n} >x }=0.
\EQN
\EEL

\prooflem{Lem:negli}  The claims follow by similar arguments as those used in the proof of Lemma \ref{Lem:2}, with $b_n^{H_0/\beta}$ replaced by $a_n$. The assumption $\beta\leq 2H-H_0$ is used to show that
\BQNY
\lim_{n\to\IF} b_n^{H_0/\beta} a_n^{-1}\sup_{0\le t\le (t_0-\vn_0) }X(t)<\IF\ \ \textrm{and} \
\lim_{n\to\IF}a_n  ^{-1}\sup_{t\geq0} (X(t)-\hat c t^\beta)  =0,\ \textrm{a.s.}.
\EQNY
The details are thus omitted.
\QED

\begin{remark}\label{Rem:XX}
It is easy to check that  the claims in \eqref{negpart-ii-1} and \eqref{negpart-ii-2} are still valid if we  remove $\sigma_0X(t)$ from the numerators and without assuming $\beta\le 2H-H_0.$ This observation is useful for the following result.
\end{remark}

\BEL\label{Lem:Gamma}
Under the assumptions of Theorem \ref{Thm:1}, we have, for any \Pe{$\vn_0 \in (0,t_0),$}
\BQNY
\overset{(k)} {\underset{i\le n}{\max} }\sup_{(t_0-\vn_0) b_n^{1/\beta}\le t\le (t_0+\vn_0) b_n^{1/\beta}}\frac{X_i (t) -c t^\beta -b_n}{ a_n   } \ \overset{d}\to \ \Lambda^{(k)},
\EQNY
as $n\to\IF$.
\EEL
\prooflem{Lem:Gamma} We need to show that, for any $x\in \R$,
\BQN\label{eq:limL}
\lim_{n\to\IF}\pk{\overset{(k)} {\underset{i\le n}{\max} }\sup_{(t_0-\vn_0) b_n^{1/\beta}\le t\le (t_0+\vn_0) b_n^{1/\beta}}\frac{X_i (t) -c t^\beta -b_n}{ a_n   }  >x} =\pk{\Lambda^{(k)} >x}.
\EQN
\Pe{First, from \eqref{eq:wMb} we have for any $x\in \R$,
\BQN\label{eq:limsupL} \nonumber
\pk{\Lambda^{(k)}>x}&=&\lim_{n\to\IF}\pk{a_n^{-1}\LT(\wM_n^{(k)} -b_n\RT)>x}\\
&\geq&\limsup_{n\to\IF}\pk{\overset{(k)} {\underset{i\le n}{\max} }\sup_{(t_0-\vn_0) b_n^{1/\beta}\le t\le (t_0+\vn_0) b_n^{1/\beta}}\frac{X_i (t) -c t^\beta -b_n}{ a_n   }  >x}.
\EQN
}
Next, \Ji{similarly to \eqref{eq:ABB} we can derive, for any $n\in\N$ and $x\in\R$, that}
\BQNY
\pk{a_n^{-1}\LT(\wM_n^{(k)} -b_n\RT)>x}& \le& \pk{\overset{(k)} {\underset{i\le n}{\max} }\sup_{(t_0-\vn_0) b_n^{1/\beta}\le t\le (t_0+\vn_0) b_n^{1/\beta}}\frac{X_i (t) -c t^\beta -b_n}{ a_n   }  >x}\\
&&+\ \pk{\max_{i\le n} \sup_{0\le t\le (t_0-\vn_0) b_n^{1/\beta}}\frac{X_i (t) -c t^\beta -b_n}{ a_n   }  >x}\\
&& +\ \pk{\max_{i\le n} \sup_{t\ge (t_0+\vn_0) b_n^{1/\beta}}\frac{X_i (t) -c t^\beta -b_n}{ a_n   }  >x},
\EQNY
\Ji{where the last two probabilities on the right-hand side  tend to 0 as $n\to\IF$, as discussed in Remark \ref{Rem:XX}}. Thus, we obtain
\BQN\label{eq:liminfL}
\liminf_{n\to\IF}\ \pk{\overset{(k)} {\underset{i\le n}{\max} }\sup_{(t_0-\vn_0) b_n^{1/\beta}\le t\le (t_0+\vn_0) b_n^{1/\beta}}\frac{X_i (t) -c t^\beta -b_n}{ a_n   }  >x} \ge \pk{\Lambda^{(k)} >x}.
\EQN
Therefore, \eqref{eq:limL} follows from \eqref{eq:limsupL} and \eqref{eq:liminfL}, and  the proof is complete.  \QED

\subsubsection{Proof for (ii) and (iii)}
First, similarly to \eqref{eq:ABB} we can derive, for any $\vn_0\in(0,t_0)$ and any $x\in\R$,
\BQNY
\pk{a_n^{-1} (M_n^{(k)} -b_n) >x} &\le& \pk{ \max_{i\le n} \sup_{0\le t\le (t_0-\vn_0) b_n^{1/\beta}}\frac{X_i (t) +\sigma_0X(t) -c t^\beta -b_n}{ a_n } >x }\\
&& +\pk{ \max_{i\le n} \sup_{ t\ge (t_0+\vn_0) b_n^{1/\beta}}\frac{X_i (t) +\sigma_0X(t) -c t^\beta -b_n}{ a_n  } >x }\\
&&+\pk{ \overset{(k)} {\underset{i\le n}{\max} }\sup_{(t_0-\vn_0) b_n^{1/\beta}\le t\le (t_0+\vn_0) b_n^{1/\beta}}\frac{X_i (t) +\sigma_0X(t) -c t^\beta -b_n}{ a_n } >x }\\
&=:&I_1(\vn_0,n,x)+I_2(\vn_0,n,x)+I_3(\vn_0,n,x).
\EQNY
From Lemma \ref{Lem:negli} we know
\BQNY
\lim_{n\to\IF} I_1(\vn_0,n,x)= \lim_{n\to\IF} I_2(\vn_0,n,x)=0.
\EQNY
For the remaining $I_3(\vn_0,n,x)$, we note that
\BQNY
I_3(\vn_0,n,x)\le  \pk{ \overset{(k)} {\underset{i\le n}{\max} }\sup_{(t_0-\vn_0) b_n^{1/\beta}\le t\le (t_0+\vn_0) b_n^{1/\beta}}\frac{X_i (t)   -c t^\beta -b_n}{ a_n } + \sup_{(t_0-\vn_0)  \le t\le (t_0+\vn_0)} \sigma_0 X (t)  \frac{ b_n^{H_0/\beta} }{ a_n   } >x }.
\EQNY
Then, by \eqref{limit-abn}, Lemma \ref{Lem:Gamma}, and the independence of the Gaussian processes $X$ and $X_i$'s, we obtain
\BQN\label{ub-scii}
\limsup_{n\to\IF} \pk{a_n^{-1} (M_n^{(k)} -b_n) >x} &\le& \lim_{\vp_0\to 0}\limsup_{n\to\IF} I_3(\vn_0,n,x)\\ \nonumber
&\leq& \pk{\Lambda^{(k)} + \frac{\sigma_0 c\beta}{H}1_{\{\beta=2H-H_0\}} \NN >x},
\EQN
where in the last inequality we used 
 that, for $\beta=2H-H_0$,
\BQNY
&&t_0^{H_0} \tau (2A^2)^{\frac{H_0-\beta}{2(\beta-H)}} =t_0^{H_0} \tau (2A^2)^{-1}\\
&&\ \  =\LT(\frac{H}{c(\beta-H)}\RT)^{H_0/\beta} \LT(\frac{2(\beta-H)}{\beta}\RT)\frac{1}{2}\LT(\frac{\beta}{\beta-H}\RT)^2 \LT(\frac{H}{c(\beta-H)}\RT)^{-2H/\beta}\\
&& \ \ =c\beta/H.
\EQNY
Next, since
\BQNY
&&\quad \pk{a_n^{-1} (M_n^{(k)} -b_n)  >x} \ge I_3(\vn_0,n,x)\\
&&\ge  \pk{ \overset{(k)} {\underset{i\le n}{\max} } \sup_{(t_0-\vn_0) b_n^{1/\beta}\le t\le (t_0+\vn_0) b_n^{1/\beta}}\frac{X_i (t) -c t^\beta -b_n}{ a_n   } -    \frac{ b_n^{H_0/\beta} }{ a_n }\sup_{(t_0-\vn_0)  \le t\le (t_0+\vn_0)} (-\sigma_0 X (t))  >x },
\EQNY
and thus by the same reason as above we have
\BQN\label{lb-scii}
\liminf_{n\to\IF} \pk{a_n^{-1} (M_n^{(k)} -b_n) >x} &\ge& \lim_{\vp_0\to 0}\liminf_{n\to\IF} I_3(\vn_0,n,x)\\ \nonumber
&=& \pk{\Lambda^{(k)} + \frac{\sigma_0c\beta}{H}1_{\{\beta=2H-H_0\}} \NN >x}.
\EQN
Consequently, combining \eqref{ub-scii} and \eqref{lb-scii} yields
\BQNY
\lim_{n\to\IF} \pk{a_n^{-1} (M_n^{(k)} -b_n)  >x} = \pk{\Lambda^{(k)} + \frac{\sigma_0 c\beta}{H}1_{\{\beta=2H-H_0\}} \NN >x}.
\EQNY
This completes the proof for both (ii) and (iii). 

\subsection{Proof of Remarks \ref{Rem:1} (d)}  We first consider the scenarios (i) and (iii) where $\beta\ge 2H-H_0$. The claim of non-existence of a phantom distribution function can be proved by a contradiction. If a phantom distribution function $G$ exists, then by definition we know for a sequence 
$u_n(x)=e_n x+b_n, n\in \N$ with $x\in \R$ (here $e_n =\sigma_0t_0^{H_0}b_n^{H_0/\beta}$ under scenario (i) and $e_n=a_n$ under scenario (iii)), 
\BQNY
\pk{M_n\le u_n(x)}-G^n(u_n(x)) \to 0, \ \ \ \ n\to\IF,
\EQNY
which, by Theorem \ref{Thm:1}, implies that
\BQNY
G^n(u_n(x)) \to 
\left\{
\begin{array}{ll}
\pk{\NN\le x}, &  \mbox{if $\beta> 2H-H_0$,}\\[0.1cm]
 \pk{\Lambda^{(1)}+\frac{\sigma_0c\beta}{H}\NN\le x},  & \mbox{if $\beta= 2H-H_0$,}
\end{array}
\right.
 \ \ \ \ n\to\IF.
\EQNY
The above result is not possible because these limiting distributions are not members of the only three  possible non-degenerate extreme value distribution families for IID sequence.  Thus, there is no phantom distribution function for the stationary sequence $\{Q_i\}_{i\ge1}$ under scenarios (i) and (iii).

The claim of existence of a continuous phantom distribution function under scenario (ii)  follows by applying Theorem 2 of \cite{DJL15}. 
Indeed, it can be shown by Theorem \ref{Thm:1}  that
\BQNY
\pk{M_{[nt]} \le b_n} =\pk{ a_{[nt]}^{-1} (M_{[nt]}-b_{[nt]}) \le a_{[nt]}^{-1} (b_n-b_{[nt]}) }   \to e^{-t}, \ \ \forall t> 0.
\EQNY

\subsection{Proof of Theorem \ref{Thm:2}} In the following two subsections, we present the proof for scenario (i) and scenarios (ii)-(iii), respectively.

\subsubsection{Proof for (i)} 
Due to the weak convergence result in scenario (i) of Theorem \ref{Thm:1} and the arguments as in the proof of Proposition 2.1 in \cite{Res87},  it is sufficient to show that
\BQN\label{eq:Llam}
\lim_{L\to\IF} \limsup_{n\to\IF} \int_L^\IF \lambda s^{\lambda-1}\pk{\LT|b_n^{-H_0/\beta} (M_n^{(k)} -b_n)\RT|>s } ds =0.
\EQN
 Note that
\BQNY
&&\quad \int_L^\IF \lambda s^{\lambda-1} \pk{\LT|b_n^{-H_0/\beta} (M_n^{(k)} -b_n)\RT|>s } ds \\
&& \Pe{=} \int_L^\IF \lambda s^{\lambda-1}  \pk{ b_n^{-H_0/\beta} (M_n^{(k)} -b_n)>s } ds + \int_L^\IF \lambda s^{\lambda-1} \pk{ b_n^{-H_0/\beta} (M_n^{(k)} -b_n) <-s } ds \\
&& \le \int_L^\IF \lambda s^{\lambda-1}  \pk{ b_n^{-H_0/\beta} (M_n^{(1)} -b_n)>s } ds + \int_L^\IF \lambda s^{\lambda-1} \pk{ b_n^{-H_0/\beta} (M_n^{(k)} -b_n) <-s } ds \\
&&=: H_1(n,L)   +   H_2(n,L).
\EQNY
Below we discuss $H_1(n,L)$ and $H_2(n,L)$ for large $n$ and $L$, and aim to find uniform (for large $n$ and large $s$) integrable upper bounds for the probability terms in their integrands so that \eqref{eq:Llam} holds.

\underline{Consider $H_1(n,L)$.} \Ji{Fix a small} $\hat c\in(0,c)$, we can choose a large enough $G$ such that
\BQN\label{eq:hatc}
\frac{(t_0(1+G))^ H}{1 + (c-\hat c) (t_0(1+G))^\beta} < \frac{t_0^H}{1+c t_0^\beta}=A,
\EQN
\BQN\label{eq:delta-G}
\delta_G:=2\LT( (1+G)^\beta (c-\hat c) /c -1\RT)>0,
\EQN
and
\BQN\label{eq:G-lowb}
\LT(\frac{c}{c-\hat c}\RT)^{-2H/\beta}   (1+\delta_G/2)^{2(1-H/\beta)}
=\LT(\frac{c}{c-\hat c}\RT)^{\Ji{-2}}(1+G)^{2(\beta-H)}  \ge 4.
\EQN
\COM{
\BQN\label{eq:three}
(1+\delta/2) \frac{c}{c-\hat c} \leq (1+G)^{\beta}
\EQN
\BQN\label{eq:delta}
(1+\delta/2)^{1-H/\beta} \sqrt{C_{2(1-H/\beta)}} > \LT(\frac{c}{c-\hat c}\RT)^{H/\beta},
\EQN
where
\BQNY
C_{2(1-H/\beta)}
=
\left\{
\begin{array}{ll}
1 & \mbox{if $\beta\geq 2H$,}\\[0.5cm]
2^{1-2H/\beta} & \mbox{if $H<\beta< 2H$.}
\end{array}
\right.
\EQNY
One can easily verify \eqref{eq:hatc}-\eqref{eq:delta}, for example, $\hat c=c/2$ and large enough $G$.
}

It follows that
\BQNY
\pk{ b_n^{-H_0/\beta} (M_n^{(1)} -b_n)  >s } &\le& \pk{ \max_{i\le n} \sup_{0\le t\le (1+G)t_0 b_n^{1/\beta}}\frac{X_i (t) +\sigma_0X(t) -c t^\beta -b_n}{ b_n  ^{H_0/\beta}} >s }\nonumber\\
&& +\pk{ \max_{i\le n} \sup_{ t\ge (1+G)t_0  b_n^{1/\beta}}\frac{X_i (t) +\sigma_0X(t) -c t^\beta -b_n}{ b_n  ^{H_0/\beta}} >s }\nonumber \\
&=:&I_{11}(n,s) +I_{12}(n,s).
\EQNY
By self-similarity, we have
\BQN\label{eq:I11}
I_{11}(n,s) &\le& \pk{ \max_{i\le n} \sup_{0\le t\le (1+G)t_0 b_n^{1/\beta}}\frac{X_i (t)  -c t^\beta -b_n}{ b_n  ^{H_0/\beta}} >s/2 } +   \pk{ \sup_{0\le t\le (1+G)t_0 b_n^{1/\beta}}\frac{\sigma_0X(t) }{ b_n  ^{H_0/\beta}} >s/2 }\nonumber\\
&\le& n \pk{\sup_{0\le t\le (1+G)t_0 }\frac{X_i (t) }{1+ct^\beta} >f(n,s)} +   \pk{ \sup_{0\le t\le (1+G)t_0  } \sigma_0X(t)   >\frac{s}{2} }, 
\EQN
where
$$
f(n,s) := \frac{s b_n^{(H_0-H)/\beta}}{2(1+c(1+G)^\beta t_0^\beta)} + b_n  ^{1-H/\beta},\ \ \ \ n\in\N, s\ge L.
$$
We have, from Proposition \ref{prop1},  that for all large $n$ and $s$
\BQN
n \pk{\sup_{0\le t\le (1+G)t_0 }\frac{X_i (t) }{1+ct^\beta} >f(n,s)} &\le& 2  \frac{A^{\frac{3}{2}-\frac{2}{\alpha}} U^{\frac{1}{\alpha}} \H_\alpha}{2^{ \frac{1}{\alpha}} B^{\frac{1}{2}}}\frac{ f(n,s)^{ -2}}{ \overset{\leftarrow}K ({f(n,s) ^{-1}}) }\exp\LT(-\LT(\frac{f(n,s) ^{2}}{2A^2}-\log n\RT)\RT) \nonumber\\  \label{part1-I11}
&\le &  \frac{A^{\frac{3}{2}-\frac{2}{\alpha}} U^{\frac{1}{\alpha}} \H_\alpha}{2^{ \frac{1}{\alpha}\Pe{-1}} B^{\frac{1}{2}}}f(n,s) ^{\gamma_0} \exp\LT(- \LT(\frac{f(n,s) ^{2}}{2A^2}-\log n\RT)\RT),
\EQN
with some $\gamma_0>1$ large enough, where the second inequality follows since  $ (v^{ 2}  \overset{\leftarrow}K ({v^{-1}}) )^{-1}, v>0$ is a regularly varying function at infinity.
 By   \eqref{tayex-bn} and the assumption $\beta>2H-H_0$, we have
\BQNY
\lim_{n\to\IF } \frac{b_n  ^{2(1-H/\beta)}/(2A^2) -\log n}{b_n  ^{\frac{\beta+H_0-2H}{\Ji{\beta}}}}=0,
\EQNY
\Pe{and thus}
\COM{ 
 By the definition of \eqref{eq:bnan}, we know
\BQNY
\frac{b_n  ^{2(1-H/\beta)}/(2A^2) -\log n}{b_n  ^{1-H/\beta}} =
\LT(b_n  ^{1-H/\beta}/(\sqrt{2}A) - \log n\RT)\frac{\LT(b_n  ^{1-H/\beta}/(\sqrt{2}A) + \log n\RT)}{b_n  ^{1-H/\beta}} = o(b_n  ^{1-H/\beta}),
\EQNY
and thus for any $s\geq L$ and large enough $n$,}
\BQN\label{lb-fns}
 \frac{f(n,s) ^{2}}{2A^2} - \log n  &\ge&  \LT(\frac{b_n  ^{2(1-H/\beta)}}{2A^2} - \log n \RT) +
 \frac{s b_n  ^{\frac{\beta+H_0-2H}{\Ji{\beta}}}}{2A^2 (1+c(1+G)^\beta t_0^\beta)}\\  \nonumber
&\ge&    \frac{L b_n  ^{\frac{\beta+H_0-2H}{\Ji{\beta}}}}{4A^2 (1+c(1+G)^\beta t_0^\beta)} + \frac{s-L }{2A^2 (1+c(1+G)^\beta t_0^\beta)}
\EQN
holds for all large $s$ and large  $n$.
Using the $C_r$ inequality (see Lemma A in Appendix) we know
\BQN\label{eq:fnsg}
f(n,s) ^{\gamma_0} \le \frac{s^{\gamma_0} b_n^{(H_0-H)\gamma_0/\beta}}{2^{\gamma_0}(1+c(1+G)^\beta t_0^\beta)^{\gamma_0}} + b_n  ^{(1-H/\beta)\gamma_0}.
\EQN
Then, substituting \eqref{lb-fns} and \eqref{eq:fnsg} into \eqref{part1-I11}, we obtain that,  for all large enough $n,$
\BQN\label{up-part1-I11}
n \pk{\sup_{0\le t\le (1+G)t_0 }\frac{X_i (t) }{1+ct^\beta} >f(n,s)}  \le  s^{\gamma_0} \exp\LT( -\frac{s }{2A^2 (1+c(1+G)^\beta t_0^\beta)}\RT).
\EQN



Furthermore, in the light of the Borell-TIS inequality, we get
\BQN \label{eq:I11:2}
 \int_L^\IF \lambda s^{\lambda-1} \pk{ \sup_{0\le t\le (1+G)t_0  } \sigma_0X(t)   >\frac{s}{2} } ds\le  \int_L^\IF \lambda s^{\lambda-1}\exp\LT(-\frac{(s/(2\sigma_0) - K_3)^2}{2(1+G)^{2H}t_0^{2H} }\RT)ds\to 0
\EQN
as $L\to\IF$, where $K_3 := \E{ \sup_{0\le t\le (1+G)t_0  } X(t) }<\IF$.  Consequently, 
it follows from \eqref{eq:I11} and \eqref{up-part1-I11}-\eqref{eq:I11:2} that
\BQNY
\lim_{L\to\IF} \limsup_{n\to\IF} \int_L^\IF \lambda s^{\lambda-1} I_{11}(n,s)ds =0.
\EQNY
Next, we show
\BQN\label{eq:I12}
\lim_{L\to\IF} \limsup_{n\to\IF} \int_L^\IF \lambda s^{\lambda-1} I_{12}(n,s)ds =0,
\EQN
which, together with the above equation, will give the desired result that
\BQN \label{eq:I1}
\lim_{L\to\IF} \limsup_{n\to\IF} H_{1}(n,L) =0.
\EQN
Now, we focus on  $I_{12}(n,s)$. 
By self-similarity, we have, for large $n$,
\BQN\label{eq:I12_2}
I_{12}(n,s)&\le&\pk{ \max_{i\le n} \sup_{ t\ge (1+G)t_0  b_n^{1/\beta}}\frac{X_i (t)   -(c-\hat c) t^\beta -b_n}{ b_n  ^{H_0/\beta}} >s/2 } \nonumber\\
&&+\ \pk{   \sup_{ t\ge (1+G)t_0  b_n^{1/\beta}}\frac{\sigma_0X (t)  -\hat c  t^\beta  }{ b_n  ^{H_0/\beta}} >s/2 } \nonumber\\ \nonumber
&\le& n \pk{  \sup_{ t\ge (1+G)t_0 } X_1 (t)  -(1+(c-\hat c) t^\beta)  b_n  ^{1-H/\beta} >b_n^{(H_0-H)/\beta} s/2  } \\ \nonumber
&&+\ \pk{   \sup_{ t\ge 0} \sigma_0 X (t)  -\hat c  t^\beta    >b_n  ^{H_0/\beta}s/2 }\\
&\leq &  n \pk{  \sup_{ t\ge (1+G)t_0 } \frac{X_1 (t)}{  d_{n,s}+(c-\hat c) t^\beta}   >  b_n  ^{1-H/\beta}  }
+\ \pk{   \sup_{ t\ge 0} \sigma_0 X (t)  -\hat c  t^\beta    >s/2 },
\EQN
with
$$
d_{n,s} :=1+\frac{s}{2}b_n  ^{H_0/\beta-1}.
$$
From Proposition \ref{prop1}, we see
\BQN \label{eq:tail_hatc}
\lim_{L\to\IF}\limsup_{n\to\IF}  \int_L^\IF \lambda s^{\lambda-1}  \pk{   \sup_{ t\ge 0} \sigma_0X (t)  -\hat c  t^\beta    >s/2 }  ds =0.
\EQN

Furthermore, define 
\BQNY
g_{n,s}(t) :=\frac{t^ H}{d_{n,s} + (c-\hat c) t^\beta}, \ \ \ \ t\ge0.
\EQNY
By \eqref{maxpointre}, we know the unique maximum point of $g_{n,s}(t), t\ge 0$ is given by
\BQNY
t_{n,s}^* = t_0  \LT(\frac{c}{c-\hat c } d_{n,s}\RT)^{1/\beta}.
\EQNY
Recalling $\delta_G$ defined in \eqref{eq:delta-G}, we have
\BQNY
s\le\delta_G b_n^{1-H_0/\beta} \ \   \Leftrightarrow    \ \  t_{n,s}^* \le (1+G)t_0.
\EQNY
Therefore, we can divide the following integral into two parts,
\BQN\nonumber
&&\quad\int_L^\IF n\lambda s^{\lambda-1} \ \pk{  \sup_{ t\ge (1+G)t_0 } \frac{X_1 (t)}{  d_{n,s}+(c-\hat c) t^\beta}   >  b_n  ^{1-H/\beta}  }  ds \\ \nonumber
&&=\LT(\int_L^{\delta_G b_n^{1-H_0/\beta}} + \int^{\infty}_{\delta_G b_n^{1-H_0/\beta}}\RT)
n\lambda s^{\lambda-1} \pk{  \sup_{ t\ge (1+G)t_0 } \frac{X_1 (t)}{  d_{n,s}+(c-\hat c) t^\beta}   >  b_n  ^{1-H/\beta}  } ds\\ \label{eq:J1J2}
&&=:J_{11}(n,L)+ J_{12}(n,L).
\EQN
For the first integral $J_{11}(n,L)$, \Ji{since} $s\le\delta_G b_n^{1-H_0/\beta}$,  we obtain
\BQNY
\sup_{ t\ge (1+G)t_0 }  g_{n,s}(t) = g_{n,s}( (1+G)t_0) =\frac{(t_0(1+G))^ H}{d_{n,s} + (c-\hat c) (t_0(1+G))^\beta},
\EQNY
and thus by the Borell-TIS inequality,
\BQNY
n \pk{  \sup_{ t\ge (1+G)t_0 } \frac{X_1 (t)}{  d_{n,s}+(c-\hat c) t^\beta}   >  b_n  ^{1-H/\beta}  } \le n \exp\LT( - \frac{ \LT(d_{n,s} + (c-\hat c) (t_0(1+G))^\beta\RT)^2  }{2(t_0(1+G))^ {2H} } (b_n ^{1-H/\beta} -K_4)^2\RT)
\EQNY
\Pe{holds} for all large $n$ such that $b_n ^{1-H/\beta} >K_4$, where $K_4:=\E{\sup_{ t\ge (1+G)t_0 } \frac{X_1 (t)}{  1+(c-\hat c) t^\beta}  }\Pe{<\IF}$. Since
\BQNY
\LT(d_{n,s} + (c-\hat c) (t_0(1+G))^\beta\RT)^2   \ge \LT(1 + (c-\hat c) (t_0(1+G))^\beta\RT)^2  + \LT(1 + (c-\hat c) (t_0(1+G))^\beta\RT) s b_n^{ H_0/\beta-1 },
\EQNY
it follows, by \eqref{tayex-bn}, \eqref{eq:hatc} and the assumption $\beta>2H-H_0$, that, for all large $n$,
\BQNY
n \pk{  \sup_{ t\ge (1+G)t_0 } \frac{X_1 (t)}{  d_{n,s}+(c-\hat c) t^\beta}   >  b_n  ^{1-H/\beta}  } &\le& n \exp\LT( - \frac{ \LT(1+ (c-\hat c) (t_0(1+G))^\beta\RT)^2  }{2(t_0(1+G))^ {2H} } (b_n ^{1-H/\beta} -K_4)^2\RT)\\
&&\times  \exp\LT( - \frac{ \LT(1+ (c-\hat c) (t_0(1+G))^\beta\RT) s   }{4 (t_0(1+G))^ {2H} }  b_n^{ (\beta+H_0-2H)/\beta } \RT)\\
&\le &  \exp\LT( - K_0s  \RT),
\EQNY
holds, with some constant $K_0>0.$
Thus,
\BQN\label{eq:LnJ1}
\lim_{L\to\IF} \limsup_{n\to\IF} J_{11}(n,L) \le \lim_{L\to\IF}   \int_L^\IF \lambda s^{\lambda-1} e^{-K_0 s}ds =0.
\EQN

For the second integral $J_{12}(n,L)$, since $s\ge\delta_G b_n^{1-H_0/\beta}$, we  have
\BQNY
\sup_{ t\ge (1+G)t_0 }  g_{n,s}(t) = g_{n,s}( t_{n,s}^*) =\LT(\frac{c}{c-\hat c}\RT)^{H/\beta} d_{n,s} ^{H/\beta -1} A,
\EQNY
and thus by Borell-TIS inequality, for large enough $n,$
\BQNY
n \pk{  \sup_{ t\ge (1+G)t_0 } \frac{X_1 (t)}{  d_{n,s}+(c-\hat c) t^\beta}   >  b_n  ^{1-H/\beta}  }  \le  n \exp\LT( - \frac{1 }{2A^2} \LT(\frac{c}{c-\hat c}\RT)^{-2H/\beta} d_{n,s} ^{2(1-H/\beta )}  (b_n ^{1-H/\beta} -K_4)^2\RT).
\EQNY
Using a change of variable
$$
v=\frac{s-\delta_G b_n^{1-H_0/\beta}}{b_n^{1-H_0/\beta}},
$$
 and the  $C_r$ inequality, 
 we get
\BQNY
d_{n,s} ^{2(1-H/\beta)}&=&(1+\delta_G/2)^{2(1-H/\beta)}\LT(1+v/(2+\delta_G)\RT)^{2(1-H/\beta)}\\
&\geq& \frac{1}{2} (1+\delta_G/2)^{2(1-H/\beta)} \LT(1+ v^{2(1-H/\beta)}/(2+\delta_G)^{2(1-H/\beta)}\RT).
\EQNY
 Therefore, for all large  $n$,
\BQNY
J_{12}(n,L) &\leq & \int^{\infty}_{\delta_G b_n^{1-H_0/\beta}} n \lambda s^{\lambda-1}\exp\LT( - \frac{1 }{2A^2} \LT(\frac{c}{c-\hat c}\RT)^{-2H/\beta} d_{n,s} ^{2(1-H/\beta )}  (b_n ^{1-H/\beta} -K_4)^2\RT) ds\\
&\leq& \lambda n b_n^{\lambda(1-H_0/\beta)} \exp\LT(- \frac{1 }{4A^2} \LT(\frac{c}{c-\hat c}\RT)^{-2H/\beta}  (1+\delta_G/2)^{2(1-H/\beta)}  (b_n ^{1-H/\beta} -K_4)^2\RT)\\
&&\quad \times \int_0^{\infty} (v+\delta_G)^{\lambda-1}  \exp\LT(- \frac{1 }{4A^2} \LT(\frac{c}{c-\hat c}\RT)^{-2H/\beta} \LT(\frac{v}{2}\RT)^{2(1-H/\beta)} (b_n ^{1-H/\beta} -K_4)^2\RT) dv\\
& \le &    \lambda n b_n^{\lambda(1-H_0/\beta)} \exp\LT(- \frac{1 }{A^2}   (b_n ^{1-H/\beta} -K_4)^2\RT)      \times \int_0^{\infty} (v+\delta_G)^{\lambda-1}  \exp\LT(-  v^{2(1-H/\beta)}  \RT) dv,
\EQNY
where in the last inequality we have used \eqref{eq:G-lowb}.
 This, together with \eqref{tayex-bn}, implies 
\BQN\label{eq:J2nl}
\lim_{L\to\IF} \limsup_{n\to\IF} J_{12}(n,L)   =0.
\EQN
Consequently, substituting \eqref{eq:LnJ1} and \eqref{eq:J2nl} into \eqref{eq:J1J2}, and recalling \eqref{eq:I12_2}-\eqref{eq:tail_hatc}, we prove the  claim in \eqref{eq:I12}. This gives us the desired result presented in \eqref{eq:I1}.

\underline{Now consider $H_2(n,L)$.} We have, by self-similarity and symmetry of Normal distribution, that
\BQNY\nonumber
\pk{ b_n^{-H_0/\beta} (M_n^{(k)} -b_n)  <-s } &\le& \pk{\overset{(k)} {\underset{i\le n}{\max} } \frac{X_i(t_0 b_n^{1/\beta}) +\sigma_0X(t_0 b_n^{1/\beta}) -(1+ct_0^\beta) b_n}{b_n^{H_0/\beta}} <-s}\\ \nonumber
&=&\pk{\sigma_0X(t_0 )+\LT(\overset{(k)} {\underset{i\le n}{\max} }\  X_i(t_0 )  -(1+ct_0^\beta) b_n^{1-H/\beta}\RT) b_n^{(H-H_0)/\beta} <-s}\\ \nonumber
&\le&\pk{X(1) <- \frac{s}{2t_0  ^{H_0}\sigma_0}} + \pk{\LT(\overset{(k)} {\underset{i\le n}{\max} }\ X_i(1 )  - \frac{1+ct_0^\beta}{t_0^H} b_n^{1-H/\beta}\RT) b_n^{(H-H_0)/\beta}<-\frac{s}{2t_0  ^H}}\\ \label{ub-H2nL}
&=:& I_{21}(s) +I_{22}(n,s).
\EQNY
Obviously,
\BQN\label{ub-I21}
\lim_{L\to\IF} \lim_{n\to\IF} \int_L^\IF \lambda s^{\lambda-1} I_{21}(s) ds \le \lim_{L\to\IF}  \int_L^\IF \lambda s^{\lambda-1}\frac{2t_0^{H_0}\sigma_0}{\sqrt{2\pi} s} \exp\LT(-\frac{s^2}{ 8 t_0^{2H_0}\sigma_0^2}\RT) ds =0.
\EQN

Recalling $Z_n^{(k)}$ and $r_n$ as defined in \eqref{eq:NLam}-\eqref{nor-asy}, we obtain, by \eqref{lim-logn-bn}
and \eqref{asy-bnH-H0}, that 
\BQNY
I_{22}(n,s) &=& \pk{\frac{Z_n^{(k)}}{\sqrt{2 \log n} b_n^{(H_0-H)/\beta}} -r_n <-\frac{s}{2t_0  ^H}}\\
&=&\pk{Z_n^{(k)} < - \sqrt{2 \log n} b_n^{(H_0-H)/\beta}\LT(\frac{s}{2t_0  ^H} -r_n\RT)}\\
&\le& \pk{Z_n^{(k)} <- s } \le \pk{\abs{ Z_n^{(k)}} >s }\\
&\le& s^{-\kappa} \E{\abs{ Z_n^{(k)}}^\kappa}
\EQNY
holds for any $\kappa>0$ and  all large $n$ and $L$,  where the last inequality follows from Markov inequality. 
Choosing $\kappa>\lambda$  and then by Proposition \ref{Lem-mom-Ynk}, we conclude that
\BQN\label{ub-I22n}
\lim_{L\to\IF} \limsup_{n\to\IF} \int_L^\IF \lambda s^{\lambda-1} I_{22}(n,s)  ds
\leq \lim_{L\to\IF} \frac{\lambda}{\kappa-\lambda}\E{\abs{\Lambda^{(k)}}^{\kappa}}  L^{\lambda-\kappa}=0.
\EQN
Therefore, combining \eqref{ub-H2nL}-\eqref{ub-I22n} yields
\BQNY
\lim_{L\to\IF} \limsup_{n\to\IF}  H_{2}(n, L)   =0,
\EQNY
which together with \eqref{eq:I1} establishes \eqref{eq:Llam}. This completes the proof for scenario (i). 

\subsubsection{Proof for (ii) and (iii)}


The idea of proof for these two scenarios is similar to that for scenario (i),  thus we shall highlight the differences and omit some of the details when similar arguments in the proof for scenario (i) are applicable here.
%
\COM{We need to prove
\BQN\label{eq:Mbn2}
\lim_{n\to\IF}\E{\LT|\frac{M_n -b_n}{a_n}\RT|^\lambda}\ = \left\{
  \begin{array}{ll}
 \E{| \Lambda|^\lambda} , & \hbox{for scenario\ } (ii) \\
\E{ \abs{  \Lambda + \frac{c\beta}{H} \NN}^\lambda}, & \hbox{for scenario\ } (iii).
  \end{array}
\right.
\EQN
Due to the weak convergence result in scenario (i) of Theorem \ref{Thm:1}, we have by the arguments as in Proposition 2.1 in \cite{Res87} that it is sufficient to show
}
It is sufficient to show
\BQN\label{eq:Llam2}
\lim_{L\to\IF} \limsup_{n\to\IF} \int_L^\IF \lambda s^{\lambda-1}\pk{\LT|a_n^{-1}\LT(M_n^{(k)} -b_n\RT)\RT|>s } ds =0.
\EQN
 Note that
\BQNY
&&\quad \int_L^\IF \lambda s^{\lambda-1}\pk{\LT|a_n^{-1}\LT(M_n^{(k)} -b_n\RT)\RT|>s } ds \\
&& \le  \int_L^\IF \lambda s^{\lambda-1} \LT(\pk{ a_n^{-1}\LT(M_n^{(1)} -b_n\RT) >s } + \pk{a_n^{-1}\LT(M_n^{(k)} -b_n\RT)<-s }\RT) ds  \\
&&=:  H_1(n, L)   +   H_2(n,L).
\EQNY
Below we shall deal with $H_1(n,s)$ and $H_2(n,s)$, separately.

\underline{Consider $H_1(n,L)$.} 
\Ji{As before,} we can choose a large $G>0$ and some small $\hat c\in(0,c)$ such that \eqref{eq:hatc}-\eqref{eq:G-lowb} hold.

It follows that
\BQNY
\pk{  a_n^{-1}\LT(M_n^{(1)} -b_n\RT) >s }&\le& \pk{ \max_{i\le n} \sup_{0\le t\le (1+G)t_0 b_n^{1/\beta}}\frac{X_i (t) +\sigma_0X(t) -c t^\beta -b_n}{ a_n} >s }\nonumber\\
&& +\ \pk{ \max_{i\le n} \sup_{ t\ge (1+G)t_0  b_n^{1/\beta}}\frac{X_i (t) +\sigma_0X(t) -c t^\beta -b_n}{ a_n   } >s }\nonumber \\
&=:&I_{11}(n,s) +I_{12}(n,s).
\EQNY
By self-similarity, we have
\BQN\label{eq:I11-2}
I_{11}(n,s) &\le& \pk{ \max_{i\le n} \sup_{t\ge 0}\frac{X_i (t)  -c t^\beta -b_n}{ a_n  } >s/2 } +   \pk{ \sup_{0\le t\le (1+G)t_0 b_n^{1/\beta}}\frac{\sigma_0X(t) }{ a_n  } >s/2 }\nonumber\\
&=& \pk{ a_n^{-1} \LT(\wM_n^{(1)} -b_n\RT) >s/2 } +   \pk{ \sup_{0\le t\le (1+G)t_0  } \sigma_0X(t)   > {a_n}{b_n^{-H_0/\beta}}s/2 }. 
\EQN
\COM{For the first term, we have, by Proposition \ref{prop1}, that
\BQNY
\pk{  a_n^{-1} \LT(\wM_n^{(1)} -b_n\RT)>s/2 }\le 2 n R(b_n+a_n s) \exp\LT(-\frac{(b_n+a_n s)^{\tau}}{2A^2}\RT)
\EQNY
with $\tau = 2-2H/\beta$, holds for all large enough $n$ and all $s>L$.
} 
For the first term, we have from the Markov inequality and \eqref{eq:wMb2} with $m_n=n$ (choosing $\kappa>\lambda$) that
\BQN\nonumber
&&\lim_{L\to\IF} \limsup_{n\to\IF} \int_L^\IF \lambda s^{\lambda-1} \pk{   a_n^{-1} \LT(\wM_n^{(1)} -b_n\RT) >s/2 } ds \\ \label{eq:Mb-kappa}
&\le& \lim_{L\to\IF} \limsup_{n\to\IF}  \frac{2^\kappa\lambda L^{\lambda-\kappa}}{\kappa-\lambda} \E{\LT| a_n^{-1} \LT(\wM_n^{(1)} -b_n\RT)\RT|^\kappa} =0.
\EQN
Next, recalling \eqref{limit-abn} and using a similar argument as in \eqref{eq:I11:2}, we obtain
\BQN \label{eq:I11:2-2}
&&\lim_{L\to\IF} \limsup_{n\to\IF} \int_L^\IF \lambda s^{\lambda-1}\pk{ \sup_{0\le t\le (1+G)t_0  } \sigma_0X(t)   > {a_n}{b_n^{-H_0/\beta}}s/2 }ds  \nonumber\\
&\le & \lim_{L\to\IF}\int_L^\IF \lambda s^{\lambda-1} \pk{ \sup_{0\le t\le (1+G)t_0  } \sigma_0X(t)   > s A^2 /(2\tau)  }ds =0.
\EQN
Consequently, by \eqref{eq:I11-2}-\eqref{eq:I11:2-2} we have
\BQNY
\lim_{L\to\IF} \limsup_{n\to\IF} \int_L^\IF \lambda s^{\lambda-1} I_{11}(n,s)ds =0.
\EQNY
In order to obtain the desired result that
\BQN \label{eq:I1-2}
\lim_{L\to\IF} \limsup_{n\to\IF} H_{1}(n,L)  =0,
\EQN
it remains to  show
\BQNY\label{eq:I12-2}
\lim_{L\to\IF} \limsup_{n\to\IF} \int_L^\IF \lambda s^{\lambda-1} I_{12}(n,s)ds =0,
\EQNY
which results from drawing the same arguments as in the proof of scenario (i), by replacing $b_n^{H_0/\beta}$ with $a_n$ and noting that $\lim_{n\to\IF} a_n =\IF$ under the assumption of scenarios  (ii) and (iii).
The details are omitted.


 \COM{
 The following coloured can be removed in the final version ......
 \\
Recalling the $\hat c$ chosen at the beginning of the proof, we have, for all $n$ large,
\BQN\label{eq:I12_2-2}
I_{12}(n,s)&\le&\pk{ \max_{i\le n} \sup_{ t\ge (1+G)t_0  b_n^{1/\beta}}\frac{X_i (t)   -(c-\hat c) t^\beta -b_n}{ a_n } >s/2 } \nonumber\\
&&+ \pk{   \sup_{ t\ge (1+G)t_0  b_n^{1/\beta}}\frac{X (t)  -\hat c  t^\beta  }{ a_n   } >s/2 } \nonumber\\
&\le& n \pk{  \sup_{ t\ge (1+G)t_0 } X_1 (t)  -(1+(c-\hat c) t^\beta)  b_n  ^{1-H/\beta} >a_nb_n^{-H/\beta} s/2  } \\
&&+ \pk{   \sup_{ t\ge 0} X (t)  -\hat c  t^\beta    >s/2 }.\nonumber
\EQN
The second term on the right-hand side can be controlled by using Proposition \ref{prop1} for all $s$ large enough, which yields that
\BQN \label{eq:tail_hatc-2}
\lim_{L\to\IF}   \int_L^\IF \lambda s^{\lambda-1}  \pk{   \sup_{ t\ge 0} X (t)  -\hat c  t^\beta    >s/2 }  ds =0.
\EQN
 Next, the first term on the right-hand side of \eqref{eq:I12_2-2} can be rewritten as
 \BQNY
 n \pk{  \sup_{ t\ge (1+G)t_0 } X_1 (t)  -(1+(c-\hat c) t^\beta)  b_n  ^{1-H/\beta} >a_nb_n^{ -H /\beta} s/2 } =n \pk{  \sup_{ t\ge (1+G)t_0 } \frac{X_1 (t)}{  d_{n,s}+(c-\hat c) t^\beta}   >  b_n  ^{1-H/\beta}  },
 \EQNY
where
$$
d_{n,s} :=1+\frac{s}{2}a_nb_n  ^{-1}.
$$
Furthermore, denote for large $n,s$,
\BQNY
g_{n,s}(t) =\frac{t^ H}{d_{n,s} + (c-\hat c) t^\beta}, \ \ \ \ t\ge0.
\EQNY
It follows that the unique maximum point of $g_{n,s}(t), t\ge 0$ is given by
\BQNY
t_{n,s}^* = t_0  \LT(\frac{c}{c-\hat c } d_{n,s}\RT)^{1/\beta}.
\EQNY
Define, for large $G$,
\BQN\label{eq:delta-2}
\delta_G:=2\LT( (1+G)^\beta (c-\hat c) /c -1\RT)>0.
\EQN
We consider the following integral
\BQN\label{eq:J1J2-2}
&&\quad\int_L^\IF n\lambda s^{\lambda-1} \ \pk{  \sup_{ t\ge (1+G)t_0 } \frac{X_1 (t)}{  d_{n,s}+(c-\hat c) t^\beta}   >  b_n  ^{1-H/\beta}  }  ds \nonumber\\
&&=\LT(\int_L^{\delta_G b_n a_n^{-1}} + \int^{\infty}_{\delta_G b_n a_n^{-1}}\RT)
n\lambda s^{\lambda-1} \pk{  \sup_{ t\ge (1+G)t_0 } \frac{X_1 (t)}{  d_{n,s}+(c-\hat c) t^\beta}   >  b_n  ^{1-H/\beta}  } ds\\
&&=:J_1(n,L)+ J_2(n,L).\nonumber
\EQN
By the definition of $\delta_G$ in \eqref{eq:delta}, we have
\BQNY
s\le\delta_G   b_n a_n^{-1}\ \   \Leftrightarrow    \ \  t_{n,s}^* \le (1+G)t_0.
\EQNY
Consider  $J_1(n,L)$. In this case $s\le\delta_G   b_n a_n^{-1}$, then we have
\BQNY
\sup_{ t\ge (1+G)t_0 }  g_{n,s}(t) = g_{n,s}( (1+G)t_0) =\frac{(t_0(1+G))^ H}{d_{n,s} + (c-\hat c) (t_0(1+G))^\beta}.
\EQNY
By the Borell-TIS inequality,
\BQNY
n \pk{  \sup_{ t\ge (1+G)t_0 } \frac{X_1 (t)}{  d_{n,s}+(c-\hat c) t^\beta}   >  b_n  ^{1-H/\beta}  } \le n \exp\LT( - \frac{ \LT(d_{n,s} + (c-\hat c) (t_0(1+G))^\beta\RT)^2  }{2(t_0(1+G))^ {2H} } (b_n ^{1-H/\beta} -K_4)^2\RT)
\EQNY
for all large $n$ such that $b_n ^{1-H/\beta} >K_4$, where $K_4:=\E{\sup_{ t\ge (1+G)t_0 } \frac{X_1 (t)}{  1+(c-\hat c) t^\beta}  }$. Note that
\BQNY
\LT(d_{n,s} + (c-\hat c) (t_0(1+G))^\beta\RT)^2   \ge \LT(1 + (c-\hat c) (t_0(1+G))^\beta\RT)^2  + \LT(1 + (c-\hat c) (t_0(1+G))^\beta\RT) s a_nb_n^{-1 },
\EQNY
and further
\BQNY
\lim_{n\to\IF} a_n b_n^{-1 } b_n ^{2(1-H/\beta)} = \frac{\beta A^2}{\beta-H}.
\EQNY
We have, for all large $n$ and all $s>L$,
\BQNY
n \pk{  \sup_{ t\ge (1+G)t_0 } \frac{X_1 (t)}{  d_{n,s}+(c-\hat c) t^\beta}   >  b_n  ^{1-H/\beta}  } &\le& n \exp\LT( - \frac{ \LT(1+ (c-\hat c) (t_0(1+G))^\beta\RT)^2  }{2(t_0(1+G))^ {2H} } (b_n ^{1-H/\beta} -K_4)^2\RT)\\
&&\times  \exp\LT( - \frac{ \LT(1+ (c-\hat c) (t_0(1+G))^\beta\RT) s   }{2 (t_0(1+G))^ {2H} }  a_n b_n^{-1 } (b_n ^{1-H/\beta} -K_4)^2\RT)\\
&\le &  \exp\LT( -  \frac{ \LT(1+ (c-\hat c) (t_0(1+G))^\beta\RT)  \beta A^2  }{4 (t_0(1+G))^ {2H} (\beta-H)}   s  \RT),
\EQNY
Thus,
\BQN
\lim_{L\to\IF} \lim_{n\to\IF} J_1(n,L) \le \lim_{L\to\IF}   \int_L^\IF \lambda s^{\lambda-1} e^{-\frac{ \LT(1+ (c-\hat c) (t_0(1+G))^\beta\RT)  \beta A^2  }{4 (t_0(1+G))^ {2H} (\beta-H)}  s}ds =0.
\EQN
Now consider  $J_2(n,L)$. In this case $s\ge\delta_G b_n a_n^{-1}$, then we have
\BQNY
\sup_{ t\ge (1+G)t_0 }  g_{n,s}(t) = g_{n,s}( t_{n,s}^*) =\LT(\frac{c}{c-\hat c}\RT)^{H/\beta} d_{n,s} ^{H/\beta -1} A.
\EQNY
By the Borell-TIS inequality, for large $n,$
\BQNY
n \pk{  \sup_{ t\ge (1+G)t_0 } \frac{X_1 (t)}{  d_{n,s}+(c-\hat c) t^\beta}   >  b_n  ^{1-H/\beta}  }  \le  n \exp\LT( - \frac{1 }{2A^2} \LT(\frac{c}{c-\hat c}\RT)^{-2H/\beta} d_{n,s} ^{2(1-H/\beta )}  (b_n ^{1-H/\beta} -K_4)^2\RT).
\EQNY
Using a change of variable
$$
v=\frac{s-\delta_G b_n a_n^{-1}}{b_n a_n^{-1}},
$$
 and the well-known $C_r$ inequality ($C_r=1\ \textrm{if}\ r\geq1,  2^{r-1}\ \textrm{otherwise}$), we have
\BQNY
d_{n,s} ^{2(1-H/\beta)}&=&(1+\delta_G/2)^{2(1-H/\beta)}\LT(1+v/(2+\delta_G)\RT)^{2(1-H/\beta)}\\
&\geq& C_{2(1-H/\beta)} (1+\delta_G/2)^{2(1-H/\beta)} \LT(1+ v^{2(1-H/\beta)}/(2+\delta_G)^{2(1-H/\beta)}\RT).
\EQNY
 Therefore, for all large $G$ and $n$,
\BQNY
J_2(n,L) &\leq & \int^{\infty}_{\delta b_n a_n^{-1}} n \lambda s^{\lambda-1}\exp\LT( - \frac{1 }{2A^2} \LT(\frac{c}{c-\hat c}\RT)^{-2H/\beta} d_{n,s} ^{2(1-H/\beta )}  (b_n ^{1-H/\beta} -K_4)^2\RT) ds\\
&\leq& \lambda n b_n^{\lambda} a_n^{-\lambda} \exp\LT(- \frac{1 }{2A^2} \LT(\frac{c}{c-\hat c}\RT)^{-2H/\beta} C_{2(1-H/\beta)} (1+\delta_G/2)^{2(1-H/\beta)}  (b_n ^{1-H/\beta} -K_4)^2\RT)\\
&&\quad \times \int_0^{\infty} (v+\delta_G)^{\lambda-1}  \exp\LT(- \frac{1 }{2A^2} \LT(\frac{c}{c-\hat c}\RT)^{-2H/\beta} C_{2(1-H/\beta)}   \LT(\frac{v}{2+\delta_G}\RT)^{2(1-H/\beta)} (b_n ^{1-H/\beta} -K_4)^2\RT) dv\\
& \le &    \lambda n b_n^{\lambda(1-H_0/\beta)} \exp\LT(- \frac{1 }{A^2}   (b_n ^{1-H/\beta} -K_4)^2\RT)      \times \int_0^{\infty} (v+\delta_G)^{\lambda-1}  \exp\LT(-  v^{2(1-H/\beta)}  \RT) dv,
\EQNY
where we used,  for the large $G$,  that $\LT(\frac{c}{c-\hat c}\RT)^{-2H/\beta} C_{2(1-H/\beta)} (1+\delta_G/2)^{2(1-H/\beta)}\ge 2$; see \eqref{eq:delta}. This implies, by inspecting the expression of $b_n$, that
\BQN\label{eq:J2nl-2}
\lim_{L\to\IF} \lim_{n\to\IF} J_2(n,L)   =0.
\EQN
Consequently, the claim in \eqref{eq:I12-2} follows from \eqref{eq:I12_2-2}-\eqref{eq:tail_hatc-2} and \eqref{eq:J1J2-2}-\eqref{eq:J2nl-2}.  This gives us the desired result presented in  \eqref{eq:I1-2}.
}


\underline{Now consider $H_2(n,L)$.} It is worth mentioning that we cannot get useful upper bounds by simply taking a single point $t_0 b_n^{1/\beta}$ as in scenario (i). Instead, we shall use a suitable interval around  $t_0b_n^{1/\beta}$ as follows. By self-similarity,
we have,  for any $\vn_0\in(0,t_0)$, 
\BQNY\nonumber
&& \pk{ a_n^{-1}\LT(M_n^{(k)} -b_n\RT) <-s }
\le \pk{\overset{(k)} {\underset{i\le n}{\max} } \sup_{(t_0-\vn_0)b_n^{1/\beta} \le t\le (t_0+\vn_0)b_n^{1/\beta}  } \frac{X_i(t ) +\sigma_0X(t) - ct^\beta- b_n}{a_n } <-s}\\
&\le& \pk{\overset{(k)} {\underset{i\le n}{\max} } \sup_{(t_0-\vn_0)b_n^{1/\beta} \le t\le (t_0+\vn_0)b_n^{1/\beta}  } \frac{X_i(t )  - ct^\beta- b_n}{a_n } - \sup_{(t_0-\vn_0)b_n^{1/\beta} \le t\le (t_0+\vn_0)b_n^{1/\beta}  }  \frac{-\sigma_0X(t )}{ a_n}   <-s}\\
&\le& \pk{- \sup_{(t_0-\vn_0)b_n^{1/\beta} \le t\le (t_0+\vn_0)b_n^{1/\beta}  }  \frac{-\sigma_0X(t )}{ a_n}   <-\frac{s}{2}}\\
&&  + \ \pk{\overset{(k)} {\underset{i\le n}{\max} } \sup_{(t_0-\vn_0)b_n^{1/\beta} \le t\le (t_0+\vn_0)b_n^{1/\beta}  } \frac{X_i(t )  - ct^\beta- b_n}{a_n } < -\frac{s}{2}}\\
&\le&\pk{\sup_{ t_0-\vn_0 \le t\le  t_0+\vn_0   }\sigma_0X(t )>  \frac{ a_n} {2 b_n^{H_0/\beta}} s }\\
&&  +\ \pk{ (1+c(t_0+\vn_0)^\beta)\frac{\bnbH}{a_n}  \LT(\overset{(k)} {\underset{i\le n}{\max} } \sup_{(t_0-\vn_0)  \le t\le (t_0+\vn_0) } \frac{X_i (t) }{1+c t^\beta }-b_n^{1-H/\beta} \RT)   <-\frac{s}{2}}\\
&=:& I_{21}(n,s) +I_{22}(n,s).
\EQNY
As shown in \eqref{eq:I11:2-2} we can obtain
\BQN\label{eq:II21-2}
 \lim_{L\to\IF} \limsup_{n\to\IF} \int_L^\IF \lambda s^{\lambda-1} I_{21}(n,s) ds=0.
\EQN
\COM{By \eqref{limit-abn}, we have, for both (ii) and (iii), that
\BQNY
 I_{21}(n,s)  \le \pk{\sup_{ t_0-\vn_0 \le t\le  t_0+\vn_0   }X(t )>   \frac{Ht_0^{H_0}}{4 c\beta } s }
\EQNY
for all large $n$, and thus, by the Borell-TIS inequality we have 
\BQN\label{eq:II21-2}
&&\quad \lim_{L\to\IF} \limsup_{n\to\IF} \int_L^\IF \lambda s^{\lambda-1} I_{21}(n,s) ds \nonumber\\
&& \le\lim_{L\to\IF}  \int_L^\IF \lambda s^{\lambda-1} \exp\LT(-\frac{\LT(\frac{Ht_0^{H_0}}{4 c\beta }s-\E{\sup_{ t_0-\vn_0 \le t\le  t_0+\vn_0   }X(t)} \RT) ^2}{ 2 t_0^{2H}}\RT) ds \\
&&=0.\nonumber
\EQN}
In order to analyse $I_{22}(n,s)$, we shall introduce some further notation.
Denote
\BQNY
\wY_i =\sup_{(t_0-\vn_0)  \le t\le (t_0+\vn_0) }  \frac{X_i (t)}{1+c t^\beta},\quad i=1,2,\ldots.
\EQNY
We obtain, from Proposition \ref{prop1}, that
\BQNY\label{eq:wYv}
\pk{\wY_i >v} = \wR(v) \exp\LT(-\frac{v^2}{2A^2}\RT)\oo, \ \ \ \ v\to\IF,
\EQNY
where
$$
\wR(v) = \frac{A^{\frac{3}{2}-\frac{2}{\alpha}}  \H_\alpha}{2^{ \frac{1}{\alpha}} B^{\frac{1}{2}}}\frac{ v^{ -2}}{ \overset{\leftarrow}K (v^{-1}) }, \ \ \ v>0.
$$
Define
\BQNY
\wb_n&:=& A (2 \log n)^{1/2}  + A ( 2 \log n)^{-1/2}   \log(\wR((2A^2\log n)^{1/2})), \ n\in \N, \\
\wa_n&:=&  A (2 \log n)^{-1/2}, \ n\in \N.\nonumber
\EQNY
By Proposition \ref{Prop:GamLim}, we have
\BQNY
\wa_n^{-1} \LT(\overset{(k)} {\underset{i\le n}{\max} }\ \wY_i -\wb_n\RT) \ \overset{d} \to \ \Lambda^{(k)}, \ \ \ \ n\to\IF.
\EQNY
\COM{
Note  that
\BQN 
&&(1+c(t_0+ \vn_0)^\beta)\frac{\bnbH}{a_n} \LT(\overset{(k)}\max_{i\le n} \wY_i- b_n^{1-H/\beta}   \RT)\nonumber\\
&&= (1+c(t_0+ \vn_0)^\beta) \frac{\bnbH}{a_n}  \wa_n \LT(\frac{\overset{(k)}\max_{i\le n} \wY_i- \wb_n }{\wa_n} +\frac{\wb_n -b_n^{1-H/\beta}}{\wa_n} \RT).
\EQN}
Next, it can be checked that
\BQNY
\frac{\bnbH}{a_n}  = \frac{\tau}{2A^2} (2A^2\log n)^{\frac{H}{\beta \tau}-\frac{1}{\tau}+1} (1+o(1))= \frac{\tau}{2A^2} (2A^2\log n)^{1/2} (1+o(1)), \ \ n\to\IF,
\EQNY
and thus
\BQNY
\lim_{n\to\IF} (1+c(t_0\pm \vn_0)^\beta) \frac{\bnbH}{a_n}  \wa_n = A  \frac{1+c(t_0\pm \vn_0)^\beta}{t_0^{H}} =: A_{\vn_0}.
\EQNY
Further, by using second-order Taylor expansion,   as $n\to\IF,$
\BQNY
b_n^{1-H/\beta} = A\sqrt{2\log n} \LT( 1+ \frac{1}{2} (\log n)^{-1} \log(R((2A^2\log n)^{\frac{\beta}{2(\beta-H)} }))     + O\LT( (\log n)^{-2} (\log(R((2A^2\log n)^{\frac{\beta}{2(\beta-H)} }) ))^2   \RT)       \RT) .
\EQNY
Moreover,  by definition
\BQNY
R((2A^2\log n)^{\frac{\beta}{2(\beta-H)} }) =\wR((2A^2\log n)^{1/2}) =\frac{A^{\frac{3}{2}-\frac{2}{\alpha}}  \H_\alpha}{2^{ \frac{1}{\alpha}} B^{\frac{1}{2}}}\frac{ (2A^2\log n)^{ -1}}{ \overset{\leftarrow}K ((2A^2\log n)^{-1/2}) }.
\EQNY
Hence, we derive that
\BQN \label{eq:wbb}
\lim_{n\to\IF} \wa_n^{-1} \LT(\wb_n -b_n^{1-H/\beta}\RT) =0,
\EQN
and thus for all large $n$, 
\BQNY
I_{22}(n,s) 
&=&\pk{ (1+c(t_0+ \vn_0)^\beta) \frac{\bnbH}{a_n}  \wa_n \LT(\frac{\overset{(k)}\max_{i\le n} \wY_i- \wb_n }{\wa_n} +\frac{\wb_n -b_n^{1-H/\beta}}{\wa_n} \RT) <-\frac{s}{2}}  \\
&\le& \pk{ \wa_n ^{-1} \LT(\overset{(k)} {\underset{i\le n}{\max} }\ \wY_i- \wb_n \RT)    <-\frac{s}{4 A_{\vn_0}}}\\
&\le& (4A_{\vn_0})^{\kappa} s^{-\kappa} \E{\abs{\wa_n ^{-1} \LT(\overset{(k)} {\underset{i\le n}{\max} }\ \wY_i- \wb_n \RT) }^\kappa }
\EQNY
holds for any $\kappa>0$.
In view of the definition of $\wY_1$, it follows that
$$
\pk{\wY_1 \le -x} \le \pk{\frac{X_1 (t_0)}{1+c t_0^\beta} \le -x} \le \frac{A}{\sqrt{2\pi} x} e^{-\frac{x^2}{2A^2}}, \ \ \ \forall\ x>0,
$$
fulfilling \eqref{assu-Y-leta},  and thus we conclude from Proposition \ref{Lem-mom-Ynk} that, for a chosen $\kappa>\lambda,$
\BQNY\label{eq:III22}
\lim_{L\to\IF} \limsup_{n\to\IF} \int_L^\IF \lambda s^{\lambda-1} I_{22}(n,s)  ds =0.
\EQNY
This, together with \eqref{eq:II21-2}, implies
\BQNY
\lim_{L\to\IF} \limsup_{n\to\IF} H_{2}(n, L) =0.
\EQNY
Consequently, from the above equation and  \eqref{eq:I1-2} we establish  \eqref{eq:Llam2}, and thus the proof for scenarios (ii) and  (iii) is complete.

\subsection{Proof of Theorem \ref{Thm:1-1}} The proof  follows from the same lines as the proofs of Theorem \ref{Thm:1} and Theorem \ref{Thm:2}, by applying Proposition \ref{propQQ} and utilising two types of  inequalities for some of the bounds therein. These two types of inequalities are akin to the following:
\begin{itemize}
\item A lower bound using $\overset{(k)} {\underset{i\le n}{\max} }\ X_i (t_0 b_n^{1/\beta})  \ge \overset{(k)} {\underset{i\le m_n}{\max} }\ X_i (t_0 b_n^{1/\beta})$ in \eqref{eq:t0NN}. 
\item An upper bound using $ \overset{(k)} {\underset{i\le n}{\max} }\ \sup_{0\le t\le(1-\vn_0)b_n^{1/\beta}}(X_i (t)   -c_i t^\beta) \le  \overset{(k)} {\underset{i\le n}{\max} }\ \sup_{0\le t\le(1-\vn_0)b_n^{1/\beta}}(X_i (t)   -c t^\beta)$  in \eqref{eq:UUU}. 
\end{itemize}
Thus, we omit the details. The proof is complete.  


\bigskip

{\bf Acknowledgement}:
We are thankful to an associate editor and two referees for their constructive suggestions which have significantly improved the manuscript.
The research of Xiaofan Peng is partially supported by  National Natural Science Foundation of China (11701070, 71871046).

\bibliographystyle{plain}

 \bibliography{gausbibruinABCD}

\section*{Appendix}


In this appendix, we present proofs for the propositions displayed in Section \ref{Sec:Gam}.  We also include  the $C_r$ inequalities that have been frequently used in our proofs.
\medskip

{\bf Proof of Proposition \ref{Prop:GamLim}}:
The proof  follows closely from some existing results. First, thanks to the  closure property of maximum domain of attraction of the Gumbel distribution under tail equivalence (see, e.g., Proposition 3.3.28 in \cite{EKM97}), the weak limit result for $k=1$ (i.e., the Gumbel limit theorem for the maximum) follows similarly to  Theorem 1.5.3 in \cite{leadbetter1983extremes} by noting that $\lim_{n\to\IF}n(1-F(\mu_n+\nu_n x))=e^{-x}, \forall x\in\R$. 
Secondly, for general fixed $k>1$ the result follows by 
an application of Theorem 2.2.2
in \cite{leadbetter1983extremes} where it is shown that for an IID sequence the convergence for maxima is equivalent to the convergence for order statistics. \QED

\COM{
\BEL\label{Lem-mom-Ynk}
Let $\{Y_i\}_{i=1}^{\IF}$ be a sequence of IID random variables, and $Y_n^{(k)}$ be the $k$th largest of $Y_1,\cdots,Y_n$. Suppose that the right tail of $Y_1$ is equivalent to \eqref{eq:Weibull} with $\tau>1$ and satisfies
\BQN\label{assu-Y-leta}
\pk{Y_1<-x} \lesssim  {\varrho}(x)e^{-\beta x^\gamma},\quad x\to\IF,
\EQN
for some positive constants $\beta, \gamma$ and regularly varying function $ \varrho(x)>0$.

Then, we have
\BQN
\lim_{L\to\IF} \limsup_{n\to\IF}\int_L^\IF   s^{\lambda-1} \pk{ \nu_n^{-1} \LT(Y_n^{(k)}-\mu_n\RT)  <-s } ds=0 \label{lim-int-Ynk}
\EQN
holds for any $\lambda>0$, where $\mu_n$ and $\nu_n$ are given in \eqref{def-mun-nun}.
\EEL

}

{\bf Proof of Proposition \ref{Lem-mom-Ynk}}:
By Proposition \ref{Prop:GamLim}
and the same arguments as those used in the proof of Proposition 2.1 of \cite{Res87},  we only need to show that
\BQNY\label{lim-int-Ynk}
\lim_{L\to\IF} \limsup_{n\to\IF} \int_L^\IF \lambda s^{\lambda-1} \pk{\abs{\nu_n^{-1} \LT(Y_n^{(k)}-\mu_n\RT) }>s } ds =0.
\EQNY
Further, note that
\BQN
\int_L^\IF \lambda s^{\lambda-1} \pk{\abs{\nu_n^{-1} \LT(Y_n^{(k)}-\mu_n\RT) }>s } ds  &\leq& \int_L^\IF \lambda s^{\lambda-1} \pk{ \nu_n^{-1} \LT(Y_n^{(k)}-\mu_n\RT)  <-s } ds\nonumber\\
&&+
\int_L^\IF \lambda s^{\lambda-1}  \pk{ \nu_n^{-1} \LT(Y_n^{(k)}-\mu_n\RT)   >s } ds \label{eq:II12nL}\\
&=:& I_1(n,L)   +   I_2(n,L).\nonumber
\EQN
We shall first focus on $I_1(n,L)$. It can be checked that (cf. Proposition 4.1.2 in \cite{EKM97})
\BQNY
\pk{\nu_n^{-1} \LT(Y_n^{(k)}-\mu_n\RT)<-s }
=\sum_{j=0}^{k-1} \  \left(\!\!\!
   \begin{array}{c}
     n  \\
    j   \\
   \end{array}\!\!\!
 \right) \left( \pk{Y_1 \geq \mu_n -\nu_n s} \right)^j \left(\pk{Y_1\leq \mu_n -\nu_n s}\right)^{n-j}.
\EQNY
By Stirling's approximation, we see that, to verify
\BQN\label{J1}
\lim_{L\to\IF} \limsup_{n\to\IF}  I_1(n,L)=0,
\EQN
it suffices to show, for any fixed $j=0,\ldots, k-1,$
\BQN\label{lim-int-NN}
\lim_{L\to\IF}\limsup_{n\to\IF} \int_{L}^{\IF} s^{\lambda-1} n^{j} \left( \pk{Y_1 \geq \mu_n -\nu_n s} \right)^j \left(\pk{Y_1\leq \mu_n -\nu_n s}\right)^{n-j} ds =0.
\EQN
We prove this equality by dividing the above integral into  three parts as follows (with $\omega_n=\sqrt{C^{-1} \log n}$):
\BQN\label{int-NN-Psi} \nonumber
&&\quad \int_{L}^{\IF}  s^{\lambda-1} n^{j} \left( \pk{Y_1 \geq \mu_n -\nu_n s} \right)^j \left(\pk{Y_1\leq \mu_n -\nu_n s}\right)^{n-j} ds\\ \nonumber
&&  = \left( \int_{L}^{\omega_n}  + \int_{\omega_n}^{(\mu_n+L)/\nu_n}  +    \int_{(\mu_n+L)/\nu_n}^{\IF} \right)  s^{\lambda-1} n^{j} \left( \pk{Y_1 \geq \mu_n -\nu_n s} \right)^j \left(\pk{Y_1\leq \mu_n -\nu_n s}\right)^{n-j} ds\\
&&=: Q_1(n,L) + Q_2(n,L)+ Q_3(n,L),
\EQN
which is always valid for  sufficiently large $n$.

We first consider $Q_1(n,L)$.  By the fact that $\log(1-x)\leq -x, x\in(0,1)$, we obtain
\BQNY
&& \quad n^{j} \LT(\pk{Y_1 > \mu_n -\nu_n s} \RT)^{j}   \LT(\pk{Y_1 \le \mu_n -\nu_n s} \RT)^{n-j}  \\
 &&\leq n^{j} \LT(\pk{Y_1 > \mu_n -\nu_n s} \RT)^{j}  \exp \left( -(n-j)  \pk{Y_1 > \mu_n -\nu_n s}  \right).
\EQNY
Below, we consider uniform bounds for $\pk{Y_1 > \mu_n -\nu_n s}$, for all large enough $n$ and all $s\in [L, \omega_n]$.
Note
\BQN\label{eq:mns2}
\sup_{s\in [L, \omega_n]}\nu_n s/\mu_n \rightarrow 0
\EQN
as $n\to\IF$, then by the tail asymptotics of $Y_1$ in \eqref{eq:Weibull} we have
\BQNY\label{eq:wYbn2}
\frac{1}{2} \rho(\mu_n -\nu_n s) \exp\LT(-C(\mu_n -\nu_n s)^\tau\RT) \le
\pk{Y_1 > \mu_n -\nu_n s} \le
2 \rho(\mu_n -\nu_n s) \exp\LT(-\Ji{C}(\mu_n -\nu_n s)^\tau\RT)
\EQNY
holds for all large enough $n$  and all $s\in[L, \omega_n]$,
and thus 
\BQNY
&& \quad n^{j} \LT(\pk{Y_1 > \mu_n -\nu_n s} \RT)^{j}   \LT(\pk{Y_1 \le \mu_n -\nu_n s} \RT)^{n-j}  \\
 &&\leq \LT(2 n \rho(\mu_n -\nu_n s) \RT)^j \exp\LT(-jC(\mu_n -\nu_n s)^\tau\RT) \\
  &&\quad \quad \times \exp \left( -  \frac{n-j}{2} \rho(\mu_n -\nu_n s) \exp\LT(-C(\mu_n -\nu_n s)^\tau\RT) \right).
\EQNY
Next, we derive uniform bounds for $(\mu_n -\nu_n s)^\tau$, for all $s\in [L, \omega_n]$.
It can be checked that
\BQN\label{ineq-mnt-1}
1-\tau_M x \le (1-x)^{\tau} \le 1-\tau_m x, \ \ \forall\ x\in[0,1],
\EQN
where  $\tau_M=\max(\tau,1), \tau_m=\min(\tau,1)$. The Taylor's expansion yields
\BQN\label{Taylor-munt}
\mu_n^\tau =  C^{-1}\log n\left(1 + \frac{\log\left(\rho\left((C^{-1}\log n)^{1/\tau}\right)\right)}{\log n} + O\left(\frac{\log\left(\rho\left((C^{-1}\log n)^{1/\tau}\right)\right)}{\log n}\right)^2 \right),\ n\to\IF.
\EQN
Thus, it follows from \eqref{eq:mns2}-\eqref{Taylor-munt}, that, for all large enough $n$ and all $s\in[L,\omega_n]$,
\BQN\label{eq:mnus}
(\mu_n -\nu_n s)^\tau &\leq& \mu_n^\tau - \tau_m \nu_n\mu_n^{\tau-1}s \nonumber\\
&\leq& C^{-1}\log n+ C^{-1} \log\left(\rho\left((C^{-1}\log n)^{1/\tau}\right)\right)+1-\tau_m s/(2C\tau),
\EQN
and
\BQNY
(\mu_n -\nu_n s)^\tau &\geq& \mu_n^\tau - \tau_M \nu_n\mu_n^{\tau-1}s\\
&\geq& C^{-1}\log n+ C^{-1} \log\left(\rho\left((C^{-1}\log n)^{1/\tau}\right)\right)-1-2\tau_M s/(C\tau).
\EQNY
Therefore, by \eqref{eq:mns2} and the Uniform Convergence Theorem (cf.   Theorem 1.5.2 in \cite{bingham1989regular}), we get
\BQNY
\lim_{n\to\IF} \sup_{s\in[L,\omega_n]} \abs{ \frac{\rho(\mu_n-\nu_ns)}{\rho\left((C^{-1}\log n)^{1/\tau}\right)} -1 } =0,
\EQNY
and thus for all large enough $n$ and all $s\in[L,\omega_n]$,
\BQNY
  \quad n^{j} \LT(\pk{Y_1 > \mu_n -\nu_n s} \RT)^{j}   \LT(\pk{Y_1 \le \mu_n -\nu_n s} \RT)^{n-j}   \leq
 4^j \exp\left(jC+\frac{2\tau_M js}{\tau}-\frac{1}{4}\exp\left(-C+\frac{\tau_m s}{2\tau}\right)\right),
\EQNY
implying
\BQN\label{lim-Q1nL}
\lim_{L\to\IF}\limsup_{n\to\IF}Q_1(n,L) \leq  \lim_{L\to\IF} \int_{L}^{\IF}
 4^j s^{\lambda-1}\exp\left(jC+\frac{2\tau_M js}{\tau}-\frac{1}{4}\exp\left(-C+\frac{\tau_m s}{2\tau}\right)\right) ds=0.
\EQN
Now consider $Q_2(n,L) $. 
Similarly as before,  we have, for all large enough $n$,
\BQNY
Q_2(n,L) &\leq& n^j ((\mu_n+L)/\nu_n)^{\lambda+1} \left(\pk{Y_1\leq \mu_n  -\nu_n\omega_n}\right)^{n-j}  \\
&\leq&  n^j ((\mu_n+L)/\nu_n)^{\lambda+1} \exp \left( - \frac{n-j}{2} \rho(\mu_n-\nu_n\omega_n) \exp\left( -C(\mu_n  -\nu_n\omega_n)^\tau \right) \right).
\EQNY
Furthermore, it follows from an application of the upper bound in \eqref{eq:mnus} with $s=\omega_n$  and the Uniform Convergence Theorem  that
\BQNY
\exp \left( - \frac{n-j}{2}\rho(\mu_n- \nu_n\omega_n) \exp\left( -C(\mu_n -\nu_n\omega_n)^\tau \right) \right)
\leq \exp\left(-\frac{1}{4}  \exp\left(-C+ \frac{\tau_m \omega_n}{2\tau} \right) \right).
\EQNY
Therefore,
\BQN\label{lim-Q2nL}
&&\quad\lim_{L\to\IF}\limsup_{n\to\IF}Q_2(n,L)\\ \nonumber
&& \leq  \lim_{L\to\IF}\limsup_{n\to\IF} (2\tau \log n)^{\lambda+1}  \exp\left(-\frac{1}{4}  \exp\left(-C+\frac{\tau_m \sqrt{C^{-1} \log n}}{2\tau} \right) +j\log n \right)=0.
\EQN

For $ Q_3(n,L)$, we have,  by assumption \eqref{assu-Y-leta} that, for any large  $L$,
\BQNY\label{lim-Q3nL}
Q_3(n,L) &\leq&  \int_{(\mu_n+L)/\nu_n}^{\IF} s^{\lambda-1} n^j \left(\pk{Y_1\leq \mu_n -\nu_n s}\right)^{n-j}  ds\\
&\leq& \text{Const.}\cdot \int_{(\mu_n+L)/\nu_n}^{\IF}  s^{\lambda-1} n^j (\nu_n s-\mu_n)^{-(n-j)\eta}  ds\\
&\leq&  \text{Const.}\cdot n^j\nu_n^{-\lambda} L^{-(n-j)\eta+1}\int_{1}^{\IF}
(\mu_n^{\lambda-1}+(tL)^{\lambda-1})  t^{-(n-j)\eta} dt\\
&\to& 0,
\EQNY
as $n\to\IF$, where in the third inequality we  used a change of variable $t=(\nu_n s - \mu_n)/L$ and the $C_r$ inequality. Thus,
\BQN\label{lim-Q3nL}
\lim_{L\to\IF}\limsup_{n\to\IF}Q_3(n,L) =0.
\EQN

Consequently, the claim \eqref{J1} follows by combing \eqref{lim-int-NN}-\eqref{int-NN-Psi} and  \eqref{lim-Q1nL}-\eqref{lim-Q3nL}.

Now, it remains to show
\BQN\label{eq:J2LIF}
\lim_{L\to\IF}\limsup_{n\to\IF} I_2(n,L) =0.
\EQN
To this end, we shall look for suitable upper bounds of $I_2(n,L)$ for all large enough $n,L$.
It follows that
\BQNY
I_2(n,L)&\le&\int_L^\IF \lambda s^{\lambda-1}  \pk{ \nu_n^{-1} \LT(Y_n^{(1)}-\mu_n\RT)   >s } ds \\
&\le & \lambda n  \int_L^\IF s^{\lambda-1}  \pk{ Y_1>\mu_n   +\nu_n s } ds\\
&=&\lambda n \pk{ Y_1>\mu_n }  \int_L^\IF s^{\lambda-1}  \frac{\pk{ Y_1>\mu_n   +\nu_n s }}{\pk{ Y_1>\mu_n } } ds.
\EQNY
It is easy to check that $\lim_{n\to\IF} n \pk{ Y_1>\mu_n }  =1$. We now proceed to find suitable uniform integrable  bounds for $ \pk{ Y_1>\mu_n   +\nu_n s }/ \pk{ Y_1>\mu_n }$, for all large $n,s$. By \eqref{eq:Weibull} we  have
\BQNY 
\frac{\pk{ Y_1>\mu_n   +\nu_n s }}{\pk{ Y_1>\mu_n } }\ \le  \ 2 \  \frac{\rho(\mu_n+\nu_n s)}{\rho(\mu_n)}\  e^{-C\LT((\mu_n+\nu_n s)^\tau-\mu_n^\tau\RT)}
\EQNY
holds for all large enough $n, s$. Using the Potter's bounds and the $C_r$ inequality, we can show that
\BQNY
 \frac{\rho(\mu_n+\nu_n s)}{\rho(\mu_n)} = \frac{\rho(\mu_n (1+\nu_n s/\mu_n))}{\rho(\mu_n)} \le D_0 (1+s^{2\abs{\gamma}})
\EQNY
holds for all large enough $n, s$, with some constant $D_0>0$ independent of $n,s$, where we recall that $\gamma$ is the regularly varying index of $\rho(\cdot)$. Thus, we obtain, for all large $n,L$,
\BQN \label{eq:sgamma}
I_2(n,L) \le 4 \lambda D_0 \int_L^\IF   s^{\lambda-1}  (1+s^{2\abs{\gamma}})  e^{-C\mu_n^\tau \LT((1+\nu_n s/\mu_n)^\tau-1\RT)} ds.
\EQN
In order to obtain upper bounds for the exponential term in \eqref{eq:sgamma}, we shall distinguish case $\tau\ge 1$ and case $\tau<1$. For the case $\tau\ge 1$, it is obvious that
$(1+\nu_n s/\mu_n )^\tau \ge (1+\nu_n s/\mu_n)$ and thus
\BQNY
 I_2(n,L)  \le 4 \lambda D_0 \int_L^\IF   s^{\lambda-1}  (1+s^{2\abs{\gamma}})  e^{-\frac{s}{2\tau} } ds
\EQNY
for all large $n$. 
This yields \eqref{eq:J2LIF} for $\tau\ge 1$.
For the case $\tau<1$, we first fix some large $L_0>0$ such that
\BQNY
(1+x)^\tau \ge 1+\frac{1}{2} x^\tau,\ \ \ \ \forall\ x> L_0,
\EQNY
and then, we choose some $a\in(0,\tau(1+L_0)^{\tau-1})$ such that
\BQNY
(1+x)^\tau \ge 1+ a x, \ \ \ \ \forall\ x\in[0,L_0].
\EQNY
From the above two inequalities, we can obtain that
\BQNY
(1+\nu_n s/\mu_n )^\tau \ge \left\{
\begin{array}{ll}
1+\frac{1}{2}(\nu_n s/\mu_n )^\tau, &  \mbox{\Pe{if $s>L_0\mu_n/\nu_n$},}\\[0.1cm]
1+a\nu_n s/\mu_n , & \mbox{\Pe{if $s\le L_0\mu_n/\nu_n$}.}\\[0.1cm]
\end{array}
\right.
\EQNY
Further, noting that $\lim_{n\to\IF}\nu_n=\IF$ for $\tau<1$, we derive, for all large $n,L$,
\BQNY
I_2(n,L) &\le& 4 \lambda D_0 \int_L^{L_0\mu_n/\nu_n}   s^{\lambda-1}  (1+s^{2\abs{\gamma}})  e^{-\frac{as}{2\tau} } ds\\
&&+4 \lambda D_0 \int_{L_0\mu_n/\nu_n}^\IF   s^{\lambda-1}  (1+s^{2\abs{\gamma}})  e^{-  s ^\tau } ds.
\EQNY
This implies \eqref{eq:J2LIF} for $\tau< 1$. Therefore, \eqref{eq:J2LIF} is established for all $\tau>0,$ and thus the proof is complete.  \QED


{\bf Proof of Proposition \ref{propQQ}:}  
If  $m_n=n$ (i.e., $c=c_i, i\ge 1$) then the claim follows immediately from Propositions \ref{Prop:GamLim} and \ref{Lem-mom-Ynk} for the IID sequence $\{\wQ_i\}_{i\ge 1}$.
We now focus on the non-stationary case where $m_n<n$. To show \eqref{eq:wMb} is equivalent to show that, for any $x\in \R$,
\BQNY 
\lim_{n\to\IF}\pk{a_{m_n}^{-1}(\wM_n^{(k)} -b_{m_n} )>x} = \pk{\Lambda^{(k)}>x}.
\EQNY
 Clearly, we have
\BQNY
\pk{a_{m_n}^{-1}(\wM_n^{(k)} -b_{m_n} )>x} \ge \pk{a_{m_n}^{-1}(\wM_{m_n}^{(k)} -b_{m_n} )>x}
\EQNY
and
\BQNY
\pk{a_{m_n}^{-1}(\wM_n^{(k)} -b_{m_n} )>x} \le \pk{a_{m_n}^{-1}(\wM_{m_n}^{(k)} -b_{m_n} )>x}+\pk{\cup_{m_n< l\le n} \LT(\wQ_l> b_{m_n} +a_{m_n} x\RT)}.
\EQNY
We have already shown that $\lim_{n\to\IF}\pk{a_{m_n}^{-1}(\wM_{m_n}^{(k)} -b_{m_n} )>x} = \pk{\Lambda^{(k)}>x}.$
Next, note that $A=A(c)$ defined in \eqref{eq:AB} as a function of $c$ is strictly decreasing.  It follows from Proposition \ref{prop1} that, for any $l>m_n,$
\BQNY
\pk{\wQ_l > b_{m_n} +  a_{m_n} x}  = o\LT( \pk{\wQ_{1}> b_{m_n} +  a_{m_n} x} \RT), \ \ \ \ n\to\IF.
\EQNY
 Thus, 
\BQNY
\pk{\cup_{m_n< l\le n} \LT(\wQ_l> b_{m_n} +a_{m_n} x\RT)} \le (n-m_n)  o\LT( \pk{\wQ_{1}> b_{m_n} +  a_{m_n} x} \RT) \to 0, \ \ \ \ n\to\IF,
\EQNY
where we use the fact that 
\BQNY
(n-m_ n) \pk{\wQ_{1} > b_{m_n} +  a_{m_n} x}  \to (1-p) e^{-x}, \ \ \ \ n\to\IF.
\EQNY
Consequently, the claim in \eqref{eq:wMb} follows.
Next we show that \eqref{eq:wMb2} can be established  similarly as Proposition \ref{Lem-mom-Ynk}. In fact, considering in formula \eqref{eq:II12nL} $Y_n^{(k)}$ to be $\wM_n^{(k)}$, $\mu_n$ to be $b_{n}$, and $\nu_n$ to be $\sigma_0 t_0^{H_0}b_{n}^{H_0/\beta}$ if $\beta>2H-H_0$, and $\nu_n=a_{n}$, otherwise, respectively, we have 
\BQNY
 I_1(n,L) &\le &\int_L^\IF \lambda s^{\lambda-1} \pk{\nu_{m_n}^{-1} \LT(\wM_{m_n}^{(k)}-\mu_{m_n}\RT)  <-s } ds, \\
I_2(n,L) &\le& \int_L^\IF \lambda s^{\lambda-1} \pk{\nu_{m_n}^{-1} \LT({\underset{i\le n}{\max} }\ \sup_{t\ge 0}(X_i (t)   -c t^\beta) -\mu_{m_n}\RT) >s } ds.
\EQNY 
The rest of the proof follows the same lines of arguments as those in the proof of Proposition \ref{Lem-mom-Ynk}. This completes the proof.
\QED


\medskip
{\it  { \bf  Lemma A. [$C_r$ inequalities]} Let $a_i,i=1,2,\ldots,n,$ and $\alpha$ be positive constants, we have
\BQNY
\left\{
\begin{array}{ll}
n^{1-\alpha}\LT( \sum_{i=1}^n a_i \RT)^\alpha \leq \sum_{i=1}^n a_i^\alpha \leq \LT(\sum_{i=1}^n a_i\RT)^\alpha, & \mbox{if $\alpha> 1$,}\\[0.5cm]
\LT(\sum_{i=1}^n a_i\RT)^\alpha  \leq \sum_{i=1}^n a_i^\alpha \leq n^{1-\alpha}\LT(\sum_{i=1}^n a_i \RT)^\alpha , & \mbox{if $\alpha\leq 1$.}
\end{array}
\right.
\EQNY
In particular,
\BQNY
n^{-1} \sum_{i=1}^n a_i^\alpha \leq \LT(\sum_{i=1}^n a_i\RT)^\alpha \leq n^\alpha \sum_{i=1}^n a_i^\alpha.
\EQNY
 }

\COM{
\prooflem{Lem:Cr} If $\alpha>1$, then
\BQNY
\frac{\sum_{i=1}^n a_i^\alpha}{\LT(\sum_{i=1}^n a_i \RT)^\alpha} = \sum_{i=1}^n
\LT( \frac{a_i}{\sum_{i=1}^n a_i }\RT)^\alpha \leq 1.
\EQNY
Further, by Jessen's inequality,
\BQNY
\sum_{i=1}^n a_i^\alpha\frac{1}{n} \geq \LT(\sum_{i=1}^n a_i \frac{1}{n}\RT)^\alpha \Longrightarrow \sum_{i=1}^n a_i^\alpha \geq n^{1-\alpha}\LT( \sum_{i=1}^n a_i \RT)^\alpha.
\EQNY
Similarly, one can verify the conclusion for $\alpha\leq1$. \QED
}

\bigskip

\end{document}